\documentclass[12pt]{article}

\usepackage[affil-it]{authblk}

\usepackage[in]{fullpage}

\usepackage{hyperref}
\usepackage{paralist}
\usepackage{siunitx}
\usepackage{adjustbox}
\usepackage{amssymb, amsmath}
\usepackage{xcolor}
\usepackage{tikz}
\usepackage{flowchart}
\usepackage[numbers]{natbib}

\usetikzlibrary{calc}
\usetikzlibrary{decorations.pathmorphing,patterns}
\usetikzlibrary{calc,patterns,decorations.markings}
\usetikzlibrary{decorations.pathreplacing}
\usetikzlibrary{shapes,arrows,spy,positioning,decorations}
\RequirePackage{amsmath, amssymb, amsthm}

\definecolor{ultramarine}{RGB}{0,32,96}
\definecolor{forestgreen}{RGB}{0,102,51}
\definecolor{darkred}{RGB} {182,19,19}

\usepackage{xspace}

\newcommand{\RNum}[1]{\uppercase\expandafter{\romannumeral #1\relax}}  
\newcommand{\meth}[1]{\texttt{\textbf{#1}}}

\newcommand{\bigO}{\mathcal{O}}

\newcommand{\onemove}{\meth{1move}\xspace}
\newcommand{\oneswap}{\meth{1swap}\xspace}

\xspaceaddexceptions{())]\}}

\newcommand{\AttendedHomeDelivery}{\textit{Attended Home Delivery}\xspace}

\newcommand{\VRPTWN}{{\tt VRPTW}\xspace}
\newcommand{\VRPTW}{{\tt (VRPTW)}\xspace}

\newcommand{\TSPTWN}{{\tt TSPTW}\xspace}

\newcommand{\TSPsTW}{{\tt (TSPsTW)}\xspace}
\newcommand{\TSPsTWN}{{\tt TSPsTW}\xspace}

\newcommand{\ANS}{{\tt (ANS)}\xspace}
\newcommand{\ANSN}{{\tt ANS}\xspace}

\newcommand{\MILPs}{{\tt MILPs}\xspace}

\newcommand{\SOPN}{{\tt SOP}\xspace}
\newcommand{\SOP}{{\tt (SOP)}\xspace}

\newcommand{\AHDN}{{\tt AHD}\xspace}
\newcommand{\AHD}{{\tt (AHD)}\xspace}

\newcommand{\MILP}{{\tt (MILP)}\xspace}
\newcommand{\MILPN}{{\tt MILP}\xspace}

\newcommand{\NeighborhoodOne}{\textit{Inside}\xspace}
\newcommand{\NeighborhoodTwo}{\textit{Outside}\xspace}

\newcommand{\SimpleInsertion}{\textit{Simple Insertion}\xspace}
\newcommand{\LocalSearch}{\textit{Local Search}\xspace}

\newcommand{\TSPTWInsertion}{{\tt TSP(s)TW} {\it{Insertion}}\xspace}

\newcommand{\NP}{\ensuremath{\mathcal{NP}}}

\newcommand{\tfeas}{\meth{TFEAS($\mathcal{A}$)}\xspace}

\DeclareMathOperator*{\windowsone}{\mathcal{W}_{\text{NO}}  }
\DeclareMathOperator*{\windowstwo}{\mathcal{W}_{\text{OV}1.5}}
\DeclareMathOperator*{\windowsthree}{\mathcal{W}_{\text{OV}3} }

\begin{document}

\title{\bf Computational Approaches for Grocery Home Delivery Services}

\author[1,2]{Christian Truden
\thanks{Corresponding author: Christian Truden, Email: \href{mailto:christian.truden@aau.at}{christian.truden@aau.at}}
}
\author[3,1]{Kerstin Maier}
\author[3]{Anna Jellen}
\author[1]{Philipp Hungerl\"ander}

\affil[1]{Department of Mathematics, Universität Klagenfurt, Klagenfurt, Austria}
\affil[2]{Lakeside Labs GmbH, Klagenfurt, Austria}
\affil[3]{MANSIO Karl Popper Kolleg, Universität Klagenfurt, Klagenfurt, Austria}

\maketitle

\begin{abstract}
  Grocery home delivery services require customers to be present when their
  deliveries arrive.
  Hence, the grocery retailer and the customer must mutually
  agree on a time window during which the delivery can be guaranteed.
  This concept is referred to as the  Attended Home Delivery \AHD problem.
  The phase during which  customers place orders, usually through a web service,
  constitutes the computationally most challenging part of the logistical processes behind such services.
  The system must determine potential delivery time windows that can be offered to incoming customers
  and incrementally build the delivery schedule as new orders are placed.
  Typically, a Vehicle Routing Problem with Time Windows forms the underlying optimization problem.
  This work is concerned with the use case given by an international grocery retailer's online  shopping service.
  We present an analysis of efficient solution methods that can be employed to \AHDN services.
  We provide several heuristic approaches for tackling the steps mentioned above.
  However, the basic framework can be easily be adapted to be used for many similar Vehicle Routing applications.
  We provide a comprehensive computational study comparing several algorithmic strategies,
  combining heuristics utilizing Local Search operations and Mixed-Integer Linear Programs, tackling the booking process.
  Finally, we analyze the scalability and suitability of the approaches.
\end{abstract}
\noindent%
{\it Keywords: attended home delivery,
grocery home delivery,
vehicle routing problem with time windows}
\newpage

\maketitle

\section{Introduction}\label{sec:intro}

Due to the rapid digitalization of the retailing sector, \AttendedHomeDelivery ~\AHD
services \cite{camsa05} have increased in importance within the e-commerce sector.
\AHDN services come into play whenever a retailing company offers an online shopping service
that requires its customers to be present when their deliveries arrive.
Hence, the retailing company and the customer must mutually agree on a
delivery time slot during which the arrival  of the delivery as well as
the presence of the customer can be assured.

This work focuses on the online shopping service of one of the world's leading grocery chains.
Grocery home delivery services are prime examples of \AHDN services
as temperature-sensitive goods cannot be dropped off at the door.
Therefore, the customer must be present to receive their order.
However, there are many more applications that follow
the same principle, such as maintenance and repair services \cite{smartMeter},
on-demand mobility services \cite{parragh_survey_2008}, or
patient home health care services \cite{fikar_home_2017}.

\AHDN services come along with several benefits for the customers,
such as nonstop opening hours of the online store,
the avoidance of traveling to the brick-and-mortar stores,
almost no interruptions of the cold chain when buying groceries, and
no carrying of heavy or bulky items.
Despite the huge potential and benefits for customers,
online grocery shopping services pose several interrelated logistics and
optimization challenges to the operator.
Basically, all actions an e-grocery retailer has to perform throughout
the planning and fulfillment process can be split into four phases.
We briefly discuss these different phases and kindly refer to Cwioro et al. \cite{Cwioro2018}
for a detailed description of the overall logistics process behind \AHDN services.
The first phase is called  \textit{Tactical Planning Phase},
it appears  several months/weeks before delivery.
During this period, the fleet of delivery vehicles is defined and
drivers are assigned to them.
As next, speaking of several weeks up to days/hours before delivery comes the  \textit{Ordering Phase},
in which the grocery chain accepts orders.
Once all orders are placed, usually days/hours before delivery,
the company prepares the accepted orders for delivery during the \textit{Preparation Phase}. See \cite{vazquez-noguerol_mixed_2021} for a model that describes the order picking at a central warehouse.
Finally, during the \textit{Delivery Phase} the delivery vehicles execute the
orders according to the delivery schedule.

The interaction with customers through the online store holds several
computational challenges.
From a computational point of view, the runtime requirements for the optimization problems occurring in the \textit{Ordering Phase} are much more challenging than in the other phases.
Thus, this work is concerned with the challenges of the  \textit{Ordering Phase}.

\begin{figure}
  \centering
  \begin{adjustbox}{scale=0.7}
    %
%

\usetikzlibrary{arrows}
\pgfdeclarelayer{background}
\pgfdeclarelayer{foreground}
\pgfsetlayers{background,main,foreground}

\tikzstyle{mypredpro}=[fill=white,draw, predproc, text width=\textwidth, text centered, minimum width=0.2\textwidth,
  minimum height=0.2\textwidth]
\tikzstyle{myterminal}=[fill=white,draw, terminal, text width=\textwidth, text centered, minimum width=0.2\textwidth,
  minimum height=0.2\textwidth]
\tikzstyle{rectangle}=[fill=white,draw, process, text width=\textwidth, text centered,
  minimum height=0.5cm]
  \tikzstyle{tex}=[text centered,
  minimum height=0.5cm]
\tikzstyle{dia}=[fill=white,draw, decision, text width=0.05\textwidth, text centered,
  minimum height=0.01\textwidth]
\tikzstyle{mydecision} = [dia, text width=0.1\textwidth, minimum width=0.01\textwidth,
  minimum height=0.02\textwidth]
  \tikzstyle{predpro} = [mypredpro, text width=0.3\textwidth, minimum width=0.2\textwidth,
  minimum height=0.05\textwidth]
   \tikzstyle{term} = [myterminal, text width=0.3\textwidth, minimum width=0.2\textwidth,
  minimum height=0.05\textwidth]
  \tikzstyle{mytext} = [tex, text width=0.01\textwidth, minimum width=0.2\textwidth,
  minimum height=0.05\textwidth]
    \tikzstyle{myhelptext} = [tex, text width=0.01\textwidth, minimum width=0.2\textwidth,
  minimum height=0.05\textwidth]
\tikzstyle{phase} = [rectangle, text width=0.3\textwidth, minimum width=0.2\textwidth,
  minimum height=0.05\textwidth]
\tikzstyle{texto} = [above, text width=0.4\textwidth]
\tikzstyle{linepart} = [draw, thick, color=black!50, -latex', dashed]
\tikzstyle{line} = [draw, thick, color=black!50, -latex']
\tikzstyle{lineWithoutArrow} = [draw, thick, color=black!50]
\tikzstyle{linepartWithoutArrow} = [draw, thick, color=black!50, dashed]
\tikzstyle{ur}=[draw, text centered, minimum height=0.01em]

\newcommand{\blockdist}{1.3}
\newcommand{\edgedist}{1.5}

\newcommand{\phase}[2]{node (p#1) [phase]
  { #1\\{\footnotesize #2}}}

  \newcommand{\predpro}[2]{node (p#1) [predpro]
  { #1\\{\footnotesize #2}}}

   \newcommand{\term}[2]{node (t#1) [term]
  { #1\\{\footnotesize #2}}}

\newcommand{\mydecision}[2]{node (d#1) [mydecision]
  {  #1\\{ #2}}}

  \newcommand{\mytext}[2]{node (p#1) [mytext]
  {  #1\\{ #2}}}
  
    \newcommand{\myhelptext}[2]{node (mt#1) [myhelptext]
  {  #1\\{ #2}}}

\newcommand{\background}[5]{%
  \begin{pgfonlayer}{background}
    \path (#1.west |- #2.north)+(-2,0.5) node (a1) {};
    \path (#3.east |- #4.south)+(+2,-0.25) node (a2) {};
    \path[fill=black!10,rounded corners, draw=black!50, dashed]
      (a1) rectangle (a2);
    \path (a1.east |- a1.south)+(3.4,-0.4) node (u1)[texto]
      {\footnotesize\textbf{#5}};
  \end{pgfonlayer}}

\begin{tikzpicture}[scale=0.7,transform shape]

  \path \predpro{Tactical Planning Phase}{};

  \path (pTactical Planning Phase.south)+(0.0,-1.5) \phase{Initialize}{};
   \path (pInitialize.south)+(0.0,-1.0) \term{Web Service}{new customer};
  \path (tWeb Service.south)+(0.0,-1.0) \phase{Get TWs}{};
  \path (tWeb Service.east)+(-1.7,0.0) \mytext{}{}; 
  \path (pGet TWs.south)+(0.0,-1.0) \phase{Set TW}{};
  \path (pSet TW.west)+(-1.0,0.0) \mytext{no}{};
  \path (pSet TW.east)+(1.0,0.0) \mytext{yes}{};
   \path (pSet TW.south)+(0.0,-2) \mydecision{Insertion}{possible?};
   \path (dInsertion.east)+(2,-0.5) \myhelptext{}{}; 
  \path (dInsertion.south)+(0.0,-1.0) \phase{Improvement}{};
  \path (pImprovement.south)+(0.0,-1.0) \term{End Ordering}{};
  \path (tEnd Ordering.south)+(0.0,-1.5) \predpro{Preparation Phase}{};
  \path (pPreparation Phase.south)+(0.0,-1.5) \predpro{Delivery Phase}{};

 \path [line] (pTactical Planning Phase.south) -- node [above] {} (pInitialize) ;
  \path [line] (pInitialize.south) -- node [above] {} (tWeb Service) ;
    \path [line] (tWeb Service.south) -- node [above] {} (pGet TWs) ;
\path [linepartWithoutArrow] (dInsertion.west) -|  node [below] {} (pno) ;
\path [linepart] (pno.north) |-  node [left] {} (pGet TWs) ;
\path [linepartWithoutArrow] (dInsertion.east) -|  node [below] {} (mt) ;
\path [linepart] (pyes.north) |-  node [right] {} (p) ;
 \path [line] (pGet TWs.south) -- node [above] {} (pSet TW) ;
 \path [linepartWithoutArrow] (pImprovement.east) -|  node [below] {} (pyes) ;
  \path [line] (pSet TW.south) -- node [above] {} (dInsertion) ;
   \path [linepart] (dInsertion.south) -- node [above] {} (pImprovement) ;
  \path [linepart] (pImprovement.south) -- node [above] {} (tEnd Ordering) ;
  \path [line] (tEnd Ordering.south) -- node [above] {} (pPreparation Phase) ;
  \path [line] (pPreparation Phase.south) -- node [above] {} (pDelivery Phase) ;

  \background{pTactical Planning Phase}{pTactical Planning Phase}{pTactical Planning Phase}{pTactical Planning Phase}{1. Tactical Planning Phase}
  \background{pTactical Planning Phase}{pInitialize}{pTactical Planning Phase}{tEnd Ordering}{2. Ordering Phase}
  \background{pPreparation Phase}{pPreparation Phase}{pPreparation Phase}{pPreparation Phase}{3. Preparation Phase}
  \background{pDelivery Phase}{pDelivery Phase}{pDelivery Phase}{pDelivery Phase}{4. Delivery Phase}
\end{tikzpicture}
  \end{adjustbox}
  \caption{Depiction of the \AHDN system with focus on the \textit{Ordering Phase}.}
  \label{fig:flowchart}
\end{figure}
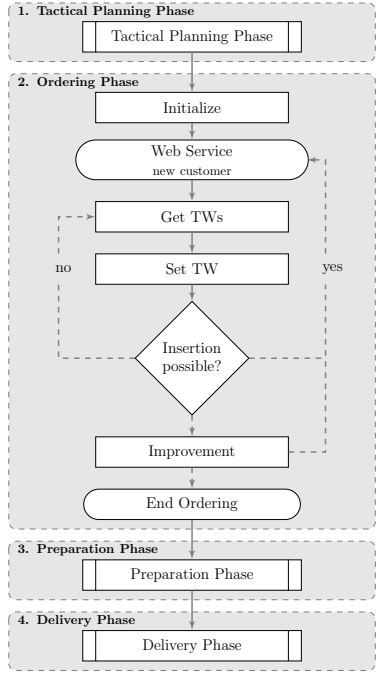

The \textit{Ordering Phase} is the phase during which customers place their grocery orders online on the company's website.
Meanwhile, the provider is challenged with solving an online variant of a \textit{Vehicle Routing Problem with Time Windows} \VRPTW,
which is known to be \NP-hard \cite{lenstra1981}.
Clearly, the website should respond to the customer requests with as little delay
as possible to ensure a smooth booking process.
Due to the tight restrictions concerning the runtime,
the naive approach of solving a new \VRPTWN instance from scratch for each new customer order
is far from being applicable in an online environment,
even when using fast Meta-heuristics.
We propose to split the  \textit{Ordering Phase}  into the following four steps (summarized in Figure \ref{fig:flowchart}).

\begin{itemize}
  \item
  \textit{Initialization step.}
  \label{step:init}
  In this step, the web service is being prepared to accept customer requests.
  Therefore, a \VRPTWN is set up, including all available vehicles with corresponding operation times.
  Since no orders have been placed yet, this results in an empty delivery schedule with a fixed fleet of vehicles.
  \item
  \textit{Get TWs step - The system determines available delivery time windows.}
  \label{step:get}
  When a new customer wants to place an order, the system has to provide all available
  delivery time windows.
  This process has to be performed in milliseconds as customers are usually
  impatient when they have to wait for technical reasons.
  Note that for calculating availabilities of time windows,
  the routing service has to calculate the travel times between all pairs of customers
  based on their provided addresses.
  The underlying mathematical problem of this step is denoted as \textit{Slot Optimization Problem} \SOP \cite{Hungerlaender2017492}.
  The \SOPN aims to determine the maximal number of available delivery time windows for a new customer.

  Optionally, for reasons of profit-maximizing, some available time windows can be
  hidden from the customer or be offered at different rates.
  In this work, we do not consider any (dynamic) slotting and pricing because
  the policy of our partnering grocery chain is to offer available time slots on a
  first-come-first-served basis, and each customer is accepted if possible.
  Nevertheless, we refer to Subsection \ref{ssec:pricing}, which provides a brief literature review about this topic.

  \item
  \textit{Set TW step - Customer books a delivery time window.}
  \label{step:set}
  Given the list of time slots (determined in the previous phase),
  the customer now chooses his or her preferred one.
  As it can take some time for the customer to decide on a time window or
  because many customers are booking simultaneously,
  the system must double-check if the selected time slot is still available.
  If the answer is yes, the customer can be added to the working schedule.
  If the answer is no, the system calls the \textit{Get TWs step} again
  to find an updated set of available time windows for the customer.
  This must be done every time a customer wants to place an order.
  Note that we do  not allow any simultaneous processing of the schedule to	avoid  queuing issues.

  \item
  \textit{Improvement step.}
  \label{step:impr}
  In the last step, optimization techniques are applied to improve the schedule.
  These are important for two reasons:
  a) offer as many time windows as possible to the customers;
  b) serve as many customers as possible.
  This can be achieved by changing the assignments of customers to vehicles and, further,
  by improving the routes of the delivery vehicles.
  As objective function we choose to minimize the total travel time as this has proven
  to be reasonable in practice.
  During very busy times, the \textit{Improvement step} can also be skipped
  or only invoked after a certain number of \textit{Set TW steps} to further improve the runtime.
  At any time of the process, we allow to have exactly
  one working schedule in the system.
\end{itemize}
We provide several heuristic approaches for tackling the steps mentioned above.
Although the approach is based on the  partnering retailing companies' requirements,
it can be easily be adapted to be used for many similar Vehicle Routing applications.

The rest of this paper is organized as follows.
In Section \ref{sec:related}, we review relevant literature.
Formal definitions of the underlying mathematical problems and descriptions of how we
conduct feasibility checks can be found in Section \ref{sec:problem}.
We propose our approaches based on \LocalSearch operations and \textit{Mixed-Integer Linear Program} \MILP formulations in
Section \ref{sec:heuristics} and
further, give suggestions on how to combine them for conducting the \textit{Get TWs},
the \textit{Set TW}, and the \textit{Improvement step}.
In Section \ref{sec:compSetup}, we present the set-up of our computational experiments
and discuss the corresponding results.
Finally, Section \ref{sec:conclude} concludes the paper.

\section{Related Work }\label{sec:related}

In this section, we briefly review relevant literature concerning
\AHDN systems as well as
solution methods for the underlying mathematical problems.

\subsection{Attended Home Delivery}

Campbell and Savelsbergh \cite{camsa05}
describe an \textit{(Attended) Home Delivery} system that decides if a new customer order is accepted.
Furthermore, the system assigns accepted orders to a time window under
consideration of the opportunity costs of the orders.
In contrast to that, in our setup the customer takes the decision
to which delivery time window her or his order is assigned to,
which requires a different setup and imposes different challenges.
Campbell and Savelsbergh describe the fulfillment process by the following three phases:
(1) order capture and promise, (2) order sourcing and assembly, and (3) order delivery.
From an algorithmic point of view, the authors propose a two-step
insertion heuristic to tackle the order capture and promise phase:
in the first step, they
employ a construction heuristic, where, starting from an empty
schedule, all already accepted requests are inserted into the schedule,
beginning with the ``heaviest'' requests.
The second step evaluates
if the new request can be inserted into the constructed
schedule in one of its acceptable time windows. Furthermore, the
authors approximate the expected profit of accepting an incoming
request. In their experimental evaluation, the heuristic provides good
results, however, only on instances with up to 100 customers, which is
much smaller than the instances we consider in our application (500 to 2000
customers).

Agatz et al.
\cite{agca08} present issues and solution approaches for the
\AHDN problem where customers select the
time window during which delivery shall take place in a very similar
setup to our problem. Their particular focus is on supermarkets that sell
groceries online. They discuss the tactical planning issues related to
the design of a time slot schedule, i.e., which time slots to offer to  customers.
Furthermore, the paper covers dynamic
time slotting (see Subsection \ref{ssec:pricing} for more details) as well as using penalties and incentives to smoothen customer demands. This work is related to our application, however, it
covers tactical considerations rather than operational difficulties
and the corresponding optimization approaches.

Han et al.
\cite{HAN2017966}
discuss an \AHDN problem that emerges
as an operational problem  at the depots of express courier companies.
The problem combines a single-depot \VRPTWN and an appointment scheduling problem.
In this setting, the couriers must arrange an appointment (via phone) to handover the delivery to the
customer. Hence,
uncertain customer behavior in responding to the arranged appointment, e.g., no-show, or random response
times, are considered.
Three main questions are addressed by the authors:
(1) The  allocation of the customers to the limited number of couriers,
(2) the sequences in which the  couriers visit the customers,
and (3) the time to meet the next customer and the maximal time to wait for the customer to show.
Consequently, the authors propose  an integrated approach that tries to balance the customers'
inconvenience and the depot's operational cost.
However, although directly related to our problem, many assumptions of this work
do not carry over to our problem setting.
In general, the preparation of grocery deliveries requires longer lead times.
Hence, the arrangement of delivery time windows takes place much earlier.
Moreover, random customer behavior is neglectable in our application.
Pan et al.
\cite{pan2017}    describe a data driven two-stage approach that focuses on predicting the
absence probability of customers for a grocery home delivery
service.
Moreover, the authors provide an excellent literature review of the online grocery
shopping process and the corresponding logistical operations.
Ehmke  \cite{Ehmke2012} gives an overview of the logistical challenges
of \AHDN systems.

As a result of a cooperation with a supermarket chain,
Vazquez-Noguerol et al. \cite{vazquez-noguerol_optimisation_2020},
present a \MILPN model for store-based e-fulfillment strategies with multiple picking locations.
The case where the orders are picked at a central warehouse is also elaborated by Vazquez-Noguerol et al. \cite{vazquez-noguerol_mixed_2021}.

\subsection{Determining Feasible Time Slots} \label{ssec:lit-GetTWs}

Most approaches in the literature \cite{Yang2016, gendreau_generalized_1998, camsa05, kohler_flexible_2019}
follow Savelsbergh's forward time slack approach \cite{savelsbergh_vehicle_1992}, which we refer to as \SimpleInsertion, for validating the feasibility of all possible delivery time windows for each incoming order.
The advantages of this approach lie in its simplicity and very short run times (see Subsection \ref{ssec:simpleinsertion} for more details).

Hungerländer et al. \cite{Hungerlaender2017492}  introduce the \SOPN.
They suggest an \textit{Adaptive Neighborhood Search} \ANS to free up time during
time windows in order to enable the insertion of new customers.
In a computational study they compare their \ANSN with two heuristics, the \SimpleInsertion
and a heuristic based on \MILPN formulations for a subproblem of the \TSPTWN,
and showed that the \ANSN is able to find much more time
slots while still fast enough for most online delivery services.
Note that their approach is restricted to non-overlapping time windows.
However, to the best of our knowledge, this is the only available paper,
which deals with a more efficient feasibility check than the \SimpleInsertion.

In this work, we adapt all mentioned heuristics such that they are able to handle
the \textit{Get TWs step} in our \AHDN considered.
As this work is not restricted to non-overlapping time windows, this results in
major adaptions for the \ANSN.

\subsection{Slotting \& Pricing}
\label{ssec:pricing}

Agatz et al. \cite{agatz_revenue_2013} discuss how proven revenue management
concepts can be translated to \AHDN services. The authors differentiate into
\textit{static} methods, i.e., forecast-based methods that are applied off-line before the actual
orders come in, and
\textit{dynamic} methods, i.e., order-based methods that are applied in real-time as
new demand comes in.
Moreover, \textit{capacity allocation} or \textit{slotting} (which time slots are made available to
which customers),
and \textit{pricing} (using delivery fees to manage customer demand),
are distinguished.
Hence, this results in the following four categories of demand management.

\textit{Differentiated (static) slotting}:
Defining the collection of delivery time windows based on
geographical regions or the preferences of customer groups.
Hence, the concentration of customer orders in a given area can be increased by
limiting the availability   of delivery options, see  \cite{Agatz2011,Hernandez2017}.

\textit{Differentiated (static) pricing}:
Differentiating between different delivery options (on a tactical level)
offered to customers by charging
different delivery fees.  Offering off-peak time discounts
or peak time premiums allows to smoothen the demand over the day, see \cite{Klein2017}.

\textit{Dynamic slotting:}
Deciding which delivery time slots to offer an incoming customer based on the
currently available capacity. More sophisticated approaches may hide delivery time
slots from unprofitable customers in order to reserve capacity for highly profitable future customers (that are predicted to arrive later on),
see \cite{cleophas_when_2014, Ehmke2014, LANG2021102305, kohler_flexible_2019}.

\textit{Dynamic pricing:}
Allows for finer levels of gradation of incentives than (dynamic) slotting.
Offering price incentives can be used to   increase the attractiveness of
time slots during which the order can be delivered more efficiently,
see \cite{asdemir_dynamic_2009, Yang2016, Klein2017b, YangStrauss2017}.

In contrast to the presented works in this subsection, we investigate the acceptance of new customer requests in terms of \textit{improving the chances of finding feasible insertions}
(given an incomplete delivery schedule) rather than
developing new acceptance criteriums for improving revenue management.

\section{Problem Description}\label{sec:problem}

In this section, we provide a detailed description of the problems we solve.
Before we state a mathematical problem description,
we recall the tasks that must  be performed in the different steps of the \textit{Ordering Phase}.

In the  \textit{Get TWs step} we aim to identify all available time windows.
The underlying optimization problem is called the Slot Optimization Problem \SOP.
Then, in the  \textit{Set TW step} we insert the customer into his or her chosen time window.
Finally, in the \textit{Improvement step} we re-optimize the schedule.
While in the  \textit{Get TWs} and the \textit{Set TW step} one must solve a feasibility version of
the \VRPTWN, the classic optimization version of the \VRPTWN must be solved in the   \textit{Improvement step}.

Thus, we start with a formal description of the \VRPTWN
in Section \ref{subsec:vrptwDescription},
continue with definitions for arrival times and feasible points of insertion
in Section \ref{ssec:arrFeas},
and  state the \SOPN finally in Section \ref{subsec:slotoptdef}.

\subsection{Vehicle Routing Problem with Time Windows} \label{subsec:vrptwDescription}
The vehicle routing problem with time windows \VRPTW is concerned with finding optimal tours for a fleet of vehicles
with given capacity constraints to deliver goods to customers within assigned
time windows.
A \VRPTWN  instance is typically defined by
a set of \textit{customers} $\mathcal{C}, ~\lvert\mathcal{C}\lvert=p$,
with corresponding \textit{order weight} function
$c \colon \mathcal{C} \to \mathbb{R}^{>0}$,
and a \textit{service time} function
$s \colon \mathcal{C} \to\mathbb R^{>0}$.
In the considered \AHDN service the individual items of an order $a \in \mathcal{C}$
are consolidated into several boxes of fixed size.
The number of required boxes  defines  the corresponding order weight $c(a)$.

Secondly, the \VRPTWN consists of a set of \textit{time windows}
$\mathcal{W} = \left\{1,\dots, q \right\}$,
where each time window  $u \in \mathcal{W}$ is
defined through the times of its begin and end $(B_u, E_u)$.
We assume that the time windows are unique, i.e.,
there do not exist time windows
$u, v \in \mathcal{W},\  u \neq v$
with $B_u = B_v$ and $E_u = E_v$.
We consider overlapping and non-overlapping time windows.
Two time windows $u$ and $v$ are \textit{non-overlapping} if and only if
$E_u \leq B_v$ or $ E_v \leq B_u$.
Further, a function $w \colon \mathcal{C} \to \mathcal{W}$ is given that
assigns to each customer a time window,
during which the delivery vehicle has to arrive at the customer.
All vehicles depart from and return to a \textit{depot} $d$.
The  \textit{travel time}  between all points $a\in V = \mathcal{C} \cup \{d\}$ is given by a function
$t \colon V  \times V \to \mathbb R^{\ge 0}$,
where we set the travel time from a customer $a$ to itself to $0$, i.e.,
$t ( a,a )=0, ~a\in V$.

Each vehicle has an assigned tour. A \textit{tour} $\mathcal{A} = ( 1,2,\dots,n)$ contains $n$ customers in the order they are visited by the vehicle and has an assigned
capacity $C_{\mathcal{A}}$.
Furthermore, each tour $\mathcal{A}$ has assigned \textit{start} and \textit{end times} that we denote as
{\it start$_{\mathcal{A}}$} and
{\it end$_{\mathcal{A}}$}, respectively.
Hence, the vehicle assigned to tour $\mathcal{A}$ can leave from the depot $d$ no earlier than
{\it start$_{\mathcal{A}}$} and must return to the depot  no later than {\it end$_{\mathcal{A}}$}.
A \textit{schedule}
$\mathcal{S}=\left\{\mathcal{A},\mathcal{B},\dots \right\}$
consists of $\lvert\mathcal{S}\rvert=m$ tours.

\subsection{Objective Function}
During the \textit{Ordering Phase} we aim to accept as many customers as possible,
while offering each customer the largest possible selection of delivery time windows.
New customer orders are accepted if at least one feasible insertion into the current delivery schedule can be found.
Ideally, the working schedule would contain few large chunks of idle time rather than many short ones.
As this is intractable to model in practice,
we alternatively choose the \textit{total travel time} as objective function
to avoid introducing an unnecessarily complicated model.
Minimizing the total travel time has proven to be a reasonable choice in practice.
Although Bent and Van Hentenryck \cite{bent_scenario-based_2004} show that the use of a consensus function in their
\textit{Multiple-Scenario Approach} results in more robust schedules and the acceptance of more customers,
their approach is not applicable to our problem, as maintaining several scenarios would introduce additional complexity and require too much computational effort.

\subsection{Arrival Times \& Feasibility}
\label{ssec:arrFeas}
For the insertion of a new customer into an existing schedule we first need to define the feasibility of an insertion.
We are going to use common concepts that we briefly describe in the following.

\subsubsection{Earliest \& Latest Arrival Times}
We consider a fixed tour $\mathcal A = (0,1, \dots,n,n+1)$,
where $0$ and $n+1$ are equal to the depot $d$ and
$(1,\ldots,n)$ are the customers assigned to the tour.
We use the concept of earliest and latest arrival time $e_i$
and $\ell_i$, as in \cite{CampSav2006},
which give the earliest (latest) time at which the
vehicle may arrive at customer $i$,
while not violating time window and travel time constraints
on the remaining tour:
\begin{align*}
  &e_0 :=start_{\mathcal A},
  ~	&&e_{j+1}:=
  \max
  \left\{
  B_{w(j+1)},~
  e_{j} + s(j)+t(j, j+1)
  \right\},
  j \in [n-1]_{0},
  \\
  &&&
  e_{n+1}:= e_n +s(n)+t(n,n+1),
  \\
  &\ell_{n+1}:=end_{\mathcal A},
  ~
  &&\ell_{j-1} :=
  \min
  \left\{
  E_{w(j-1)}, ~
  \ell_j-t(j-1,j)  - s(j-1)
  \right\},
  j \in [n] \setminus \{1\},
  \\
  &&&  \ell_0 := \ell_1 - t(0,1).
\end{align*}
Here, we assume $s(0)=s(n+1)=0$.
\\
Following the definitions,
vehicles always leave as early as possible from the depot.
This generates unnecessary idle time before serving the first customer of a tour.
Hence, once the delivery schedule is finalized,
we alter the start times of the vehicle in order to avoid this.

A schedule $\mathcal{S}$ is \textit{feasible} if all its tours are feasible.
A tour $\mathcal{A}$ is \textit{feasible} if it satisfies both of the following conditions:
\begin{align*}
  &e_i \leq E_{w(i)},  ~ i \in [n]
  \quad \wedge   \quad  e_{n+1} \leq end_{\mathcal{A}},
  \qquad&
  \text{\meth{TFEAS($\mathcal{A}$)}}, \\
  &\sum\limits_{i \in [n]} c(i)  \leq C_{\mathcal{A}},  \qquad
  &
  \text{\meth{CFEAS($\mathcal{A}$)}}.
\end{align*}
While \meth{TFEAS($\mathcal{A}$)} ensures that the arrival times at each customer assigned to tour $\mathcal{A}$ are within their assigned time windows, \meth{CFEAS($\mathcal{A}$)} guarantees that the capacity of $\mathcal{A}$ is not exceeded. Note that we do not need to check for \meth{TFEAS($\mathcal{A}$)} if $B_{w(i)} \leq  e_i,\  i \in [n],$  as this is ensured by the definition of $e_i$.

\subsubsection{Insertion Points}
The set $\Theta  (j, \mathcal{A})$ is required for the approaches applied in the \textit{Get TWs Step} and
defines after which customers we try to insert customer $j$ into $\mathcal{A}$ during its (pre)assigned time slot.
Accordingly, we define:
\begin{align*}
  &\Theta  (j, \mathcal{A}):=
  [\Theta^{-}  (j, \mathcal{A}),
  \Theta^{+}  (j, \mathcal{A})
  ],
  \intertext{where}
  &  \Theta^{-}  (j, \mathcal{A})
  :=
  \min\limits_{i \in[n]_0}
  \{
  i \colon
  B_{w(j)}
  + s(j)    \leq  \ell_{i}
  \},
  \\
  &
  \Theta^{+}  (j, \mathcal{A})
  :=
  \max\limits_{i \in[n]_0}
  \{ i  \colon
  e_{i} + s(i)   \leq E_{w(j)}
  \}    .
\end{align*}

The index  $\Theta^{-} (j, \mathcal{A})$ defines the first customer (or the depot) on tour $\mathcal{A}$ after which customer $j$ could potentially be inserted.
Likewise,   $\Theta^{+} (j, \mathcal{A})$ defines the last customer (or the depot).
Clearly, if $ \Theta^{-}  (j, \mathcal{A}) >  \Theta^{+}  (j, \mathcal{A})$,
the insertion of $j$ during $w(j)$ is infeasible.

\subsubsection{Feasibility of an Insertion}

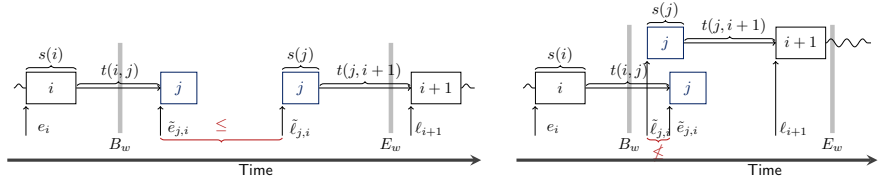
\begin{figure}[ht]
  \centering
  \scalebox{.6}{
  %
%
\tikzset{every picture/.style={line width=0.4pt}}

\begin{tikzpicture}[every node/.style={font=\sffamily\footnotesize}]

\draw[color=lightgray, line  width=1mm,-] (3cm,-1cm)  -- (3cm,1cm);
\draw[color=lightgray,  line  width=1mm,-] (9cm ,-1cm)  -- (9cm,1cm);
\node (text) at  (3cm, -1.25cm)  {$B_{w}$};
\node (text) at  (9cm, -1.25cm)  {$E_{w}$};

\node [	draw,xshift=1.47cm,yshift=0cm,	minimum width=1.1cm,
minimum height=0.7cm,] (i) {$i$};
\node [	draw,xshift=10cm,yshift=0cm,	minimum width=1.1cm,
minimum height=0.7cm,] (i1) {$i+1$};
\node [color=ultramarine,	draw,xshift=7cm,yshift=0cm,	minimum
width=0.8cm,		minimum height=0.7cm,] (j1) {$j$};
\node [color=ultramarine,	draw,xshift=4.3cm,yshift=0cm,	minimum
width=0.8cm,		minimum height=0.7cm,] (j) {$j$};

\node (inv_2) at (11.0cm,0cm) {};
\node (inv_1) at (0.5cm,0cm) {};

\draw[draw, ->] (j1.east) -> (i1.west)  ;
\draw[draw, ->] (i.east) -> (j.west) ;
\draw[decorate,decoration=snake] (inv_1.east) -> (i.west) ;
\draw[decorate,decoration=snake] (i1.east) -> (inv_2.west) ;

\draw[decorate,decoration=brace, color=black] ([shift={(0.02,0.05)}]
i.east) -- ([shift={(-0.02,0.05)}] j.west)  node[above, midway ]
{$\quad t(i,j)$ \hspace*{0.15cm} } (0,0);
\draw[decorate,decoration=brace, color=black] ([shift={(0.02,0.05)}]
j1.east) -- ([shift={(-0.02,0.05)}] i1.west)  node[above, midway ]
{$\quad t(j,i+1)$ \hspace*{0.05cm} } (0,0);
\draw[decorate,decoration=brace, color=black] ([shift={(0,0.05)}]
j1.north west) -- ([shift={(0,0.05)}] j1.north east)  node[above,
midway ] {$s(j)$} (0,0);
\draw[decorate,decoration=brace, color=black] ([shift={(0,0.05)}]
i.north west) -- ([shift={(0,0.05)}] i.north east)  node[above,
midway ] {$s(i)$} (0,0);

\draw[draw, ->] ([shift={(0,-0.7)}]j1.south  west) --
([shift={(0,-0.05)}]j1.south  west);
\node (text) at ([shift={(0.4,-0.55)}] j.south west)  {$\tilde{e}_{j,i}$};
\draw[draw, ->] ([shift={(0,-0.7)}]j.south west) --
([shift={(0,-0.05)}]j.south west);
\node (text) at ([shift={(0.4,-0.55)}]j1.south  west)  {$\tilde{\ell}_{j,i}$ };
\draw[draw, ->] ([shift={(0,-0.7)}]i.south  west) --
([shift={(0,-0.05)}]i.south  west);
\node (text) at ([shift={(0.4,-0.55)}] i.south west)  {$e_i$};
\draw[draw, ->] ([shift={(0,-0.7)}]i1.south  west) --
([shift={(0,-0.05)}]i1.south  west);
\node (text) at ([shift={(0.4,-0.55)}] i1.south west){$\ell_{i+1}$};

\draw[decorate,decoration=brace, color=darkred]
([shift={(0,-0.75)}]j1.south  west)
 --
([shift={(0,-0.75)}]j.south west)
  node[above, midway ]
{$\leq$} (0,0);

\draw[color= darkgray, draw,line  width=0.7mm, ->,>=stealth] (0.5cm,-1.6cm) -- (11cm, -1.6cm);
\draw (6cm,-1.8cm) node {Time};

\end{tikzpicture}
\begin{tikzpicture}[every node/.style={font=\sffamily\footnotesize}]

\draw[color=lightgray, line  width=1mm,-] (3cm,-1cm)  -- (3cm,1.4cm);
\draw[color=lightgray,  line  width=1mm,-] (7.5cm ,-1cm)  -- (7.5cm,1.4cm);
\node (text) at  (3cm, -1.25cm)  {$B_{w}$};
\node (text) at  (7.5cm, -1.25cm)  {$E_{w}$};

\node [	draw,xshift=1.47cm,yshift=0cm,	minimum width=1.1cm,
minimum height=0.7cm,] (i) {$i$};
\node [	draw,xshift=6.8cm,yshift=1cm,	minimum width=1.1cm,
minimum height=0.7cm,] (i1) {$i+1$};

\node [color=ultramarine,	draw,xshift=4.3cm,yshift=0cm,	minimum
width=0.8cm,		minimum height=0.7cm,] (j) {$j$};
\node [color=ultramarine,	draw,xshift=3.8cm,yshift=1cm,	minimum
width=0.8cm,		minimum height=0.7cm,] (j1) {$j$};

\node (inv_2) at (8.5cm,1cm) {};
\node (inv_1) at (0.5cm,0cm) {};

\draw[draw, ->] (j1.east) -> (i1.west)  ;
\draw[draw, ->] (i.east) -> (j.west) ;
\draw[decorate,decoration=snake] (inv_1.east) -> (i.west) ;
\draw[decorate,decoration=snake] (i1.east) -> (inv_2.west) ;

\draw[decorate,decoration=brace, color=black] ([shift={(0.02,0.05)}]
i.east) -- ([shift={(-0.02,0.05)}] j.west)  node[above, midway ]
{$\quad t(i,j)$ \hspace*{0.15cm} } (0,0);
\draw[decorate,decoration=brace, color=black] ([shift={(0.02,0.05)}]
j1.east) -- ([shift={(-0.02,0.05)}] i1.west)  node[above, midway ]
{$\quad t(j,i+1)$ \hspace*{0.05cm} } (0,0);

\draw[decorate,decoration=brace, color=black] ([shift={(0,0.05)}]
j1.north west) -- ([shift={(0,0.05)}] j1.north east)  node[above,
midway ] {$s(j)$} (0,0);
\draw[decorate,decoration=brace, color=black] ([shift={(0,0.05)}]
i.north west) -- ([shift={(0,0.05)}] i.north east)  node[above,
midway ] {$s(i)$} (0,0);

\draw[draw, ->] ([shift={(0,-0.7)}]j.south west) --
([shift={(0,-0.05)}]j.south west);
\node (text) at ([shift={(0.4,-0.55)}] j.south west)  {$\tilde{e}_{j,i}$ };
\draw[draw, ->] ([shift={(0,-1.7)}]j1.south  west) --
([shift={(0,-0.05)}]j1.south  west);
\node (text) at ([shift={(0.3,-1.55)}]j1.south  west)  {$\tilde{\ell}_{j,i}$ };
\draw[draw, ->] ([shift={(0,-0.7)}]i.south  west) --
([shift={(0,-0.05)}]i.south  west);
\node (text) at ([shift={(0.4,-0.55)}] i.south west)  {$e_i$};
\draw[draw, ->] ([shift={(0,-1.7)}]i1.south  west) --
([shift={(0,-0.05)}]i1.south  west);
\node (text) at ([shift={(0.4,-1.55)}] i1.south west) {$\ell_{i+1}$};

\draw[decorate,decoration=brace, color=darkred]
([shift={(0,-0.75)}]j.south west)
--
([shift={(0,-1.75)}]j1.south  west)
node[below, midway ] {$\nleq$} (0,0);

\draw[color= darkgray, draw,line  width=0.7mm, ->,>=stealth] (0.5cm,-1.6cm) -- (8.5cm, -1.6cm);
\draw (6cm,-1.8cm) node {Time};

\end{tikzpicture}}
  \caption{On the left \meth{TFEAS($j,i,\mathcal{A}$)} holds, on the right it does not.}
  \label{fig:feasInsert}
\end{figure}

The feasibility of the possible insertion points $\Theta  (j, \mathcal{A})$ can
be  checked easily with earliest and latest arrival times, e.g., see \cite{CampSav2006}.
Similar to the definitions above and with the help of $e_i$ and $\ell_i,\ i \in \mathcal{A}$,
we define the earliest and latest arrival time $\tilde{e}_{j,i}$ and $\tilde{\ell}_{j,i}$
for inserting a new customer $j \notin \mathcal{A}$ after $i \in \Theta  (j, \mathcal{A}) \subseteq \mathcal{A}$
within the (pre)assigned time window $w(j)$ as follows.
\begin{align*}
  &\tilde{e}_{j,i}:=
  \max
  \left\{
  B_{w(j)},~
  e_{i} + s(i)+t(i, j)
  \right\},
  \\
  &\tilde{\ell}_{j,i}:=
  \min
  \left\{
  E_{w(j)}, ~
  \ell_{i+1}-t(j,i+1)  - s(j)
  \right\}.
\end{align*}
Thus, customer $j$ can be inserted between customers $i$ and $i+1$,
such that $j$ and all subsequent customers of $\mathcal{A}$
can be served within their assigned time windows
if and only if the following condition holds.
\begin{align*}
  &\tilde{e}_{j,i} \leq \tilde{\ell}_{j,i},&&
  \meth{TFEAS($j, {i}, \mathcal{A}$)},
\end{align*}
We refer to Figure \ref{fig:feasInsert} for an illustration.

Additionally, we must check if the sum of the weights of the customer orders assigned to tour $\mathcal{A}$ does not exceed the capacity $C_\mathcal{A}$.
The insertion of $j$ into tour $\mathcal{A}$ is feasible with respect to capacity if the following condition holds:
\begin{align*}
  &\sum\limits_{i \in [n]} c(i)  + c(j) \leq C_\mathcal{A},&&
  \meth{CFEAS($j,\mathcal{A}$)}. 
\end{align*}

Assuming that all earliest (latest) arrival times and the sum of order weights on a tour $\mathcal{A}$
have already been calculated, \meth{TFEAS($j,i,\mathcal{A}$)} and \meth{CFEAS($j,\mathcal{A}$)}
allow to check the feasibility of a possible insertion 
into a given time window in $\bigO(1)$ time provided that the sequence of $\mathcal{A}$ (except for  $j$) stays the same.
If the feasibility check was successful and we decide to insert $j$,
we obtain a new tour $\mathcal{\tilde{A}}= \mathcal{A} +_i j= (0,1,\dots,i,j,i+1,\dots,n,n+1)$.
Customer  $j$ is then assigned to index $i+1$, and all indices of succeeding customers are incremented by one.
Clearly, the earliest and latest arrival times and the sum of order weights of the
modified tour must be updated.
This requires $\bigO(n)$ time \cite{CampSav2004}.

In general, there exist cases where a feasible  insertion of $j$ is only possible when the order of  $\mathcal{A}$
is changed, which makes the problem  \NP-hard.

\subsection{Slot Optimization Problem}
\label{subsec:slotoptdef}

With the introduced notation, we can recall the formal definition of the
\SOP \cite{Hungerlaender2017492},
which arises in the \textit{Get TWs step}.
We are given a feasible schedule $\mathcal{S}$ containing all scheduled customers $\mathcal{C}$,
a new customer $j,\ j \notin \mathcal{C}$, and the given set of  time slots $\mathcal{W}$.
Then, the \SOPN asks for the largest set of time slots
$\mathcal{T}_{j} \subseteq \mathcal{W}$
such that  $j$ can be serviced during each delivery slot $u \in \mathcal{T}_{j}$
by at least one vehicle of the fleet,
while assuring that all other scheduled orders stay within
their assigned time slot.
Hence, the objective is to maximize $\lvert\mathcal{T}_{j}\rvert$.

In more detail, the \SOPN aims to find at  least one feasible schedule for each of the \VRPTWN instances
consisting of  scheduled customers $\mathcal{C}$ and the new customer $j$ being temporarily assigned to
one of the time windows $u \in \mathcal{W}$.
Choosing one delivery slot for a new customer order makes the \SOPN equivalent to the
feasibility version of an appropriate \VRPTWN instance. As the
\VRPTWN is strongly \NP-hard  \cite{mako97} also the \SOPN is strongly \NP-hard
and consists of several feasibility problems that are all strongly \NP-complete.
\section{Mathematical Approaches}\label{sec:alg}
\label{sec:heuristics}

We start with introducing the algorithms that we use to tackle our proposed
\AHDN system in Sections \ref{ssec:simpleinsertion} - \ref{ssec:exact}.
Then, in Subsection \ref{subSec: solutionApproaches} we describe how they
are combined and applied to the different steps of the \textit{Ordering Phase}.

\subsection{Simple Insertion Heuristic}\label{ssec:simpleinsertion}

The \SimpleInsertion, based on \cite{savelsbergh_vehicle_1992},
takes a new customer $j$, a tour $\mathcal{A}$,
and tries to insert $j$ into the temporarily assigned time window $w(j) \in \mathcal{W}$.
It stops as soon as it finds a feasible insertion point $i \in \Theta(j, \mathcal{A})$,
i.e., when \meth{TFEAS($j,i,\mathcal{A}$)} and  \meth{CFEAS($j,\mathcal{A}$)} hold.
Note, since  the order of customers is not altered, the procedure has a linear run time
$\bigO(\lvert\mathcal{A}\rvert)$.

We iteratively apply the \SimpleInsertion to all time windows $u \in \mathcal W$
and all tours $\mathcal{A}\in\mathcal{S}_{u}$ to calculate the set of time windows that can be offered to the new customer $j$.
$\mathcal{S}_{u}$ defines the set of tours including time window $u$, i.e.,
$\mathcal{A} \in \mathcal{S}_{u}$ if and only if $start_{\mathcal A} \leq B_u < E_u \leq end_{\mathcal A}$.
A time window is considered as being \textit{available} if at least one feasible insertion point can be found.

\subsection{Local Search Heuristic}\label{sec:localSearchHeuristic}

We apply a \LocalSearch heuristic that uses the following neighborhoods for
exchanging customer orders between two tours.

\begin{compactenum}
  \item The \onemove~ neighborhood moves a customer $j$ from a tour $\mathcal{A}$
  to another tour $\mathcal{B}$, $\mathcal{A}\neq \mathcal{B}$. \\
  If at least one feasible insertion position for $j$ in $\mathcal{B}$ is found, i.e.,
  \meth{TFEAS($j,i,\mathcal{B}$)} and \meth{CFEAS($j,\mathcal{B}$)} hold,
  which additionally decreases the total travel time of the delivery schedule,
  we denote the \onemove ~ as \textit{improving}.
  \item  The \oneswap~ neighborhood exchanges two customers between two different tours, e.g.,
  switch $j \in \mathcal{A}$	with $i \in \mathcal{B}$.
  Again, if a feasible swap with decreased total travel time is found,
  we denote the \oneswap ~ as \textit{improving}.
\end{compactenum}
Savelsbergh \cite{savelsbergh_vehicle_1992} uses similar neighborhoods calling them  \textit{Relocate} and \textit{Exchange}.
Clearly, if no improving \onemove  ~(\oneswap)  could be found for a pair of
tours  $(\mathcal{A}, \mathcal{B}) \in \mathcal{S}$, $\mathcal{A} \neq \mathcal{B}$,
and both tours have not been modified meanwhile,
then there is no need to perform those operations for this pair of tours again.
Preliminary experiments showed that the computation times are reduced by a third by storing this information during the updates.

\subsection{Adaptive Neighborhood Search Heuristic}

We extend the \ANSN for solving the \SOPN
proposed by Hungerl\"ander et al. \cite{Hungerlaender2017492} such that
it can also be applied to overlapping time windows.
This results in different interdependencies between time windows as well as
slightly weaker (in)feasibility conditions.
In the following, we state all definitions required to describe our \ANSN.

\subsubsection{First/Last Customer}
For a given tour $\mathcal{A}$ we define the \textit{first} and \textit{last} customer belonging to a given time slot $u \in \mathcal{W}$ as
\begin{align*}
  &f(u):= \min_{i \in [n]}   \left\{ i: B_{u} \leq B_{w(i) }
  \leq E_{w(i)}  \leq E_{u}
  \right\},\\
  & l(u):= \max_{i \in [n]}   \left\{ i: B_{u} \leq B_{w(i) }
  \leq E_{w(i)}  \leq E_{u} \right\}.
\end{align*}
If above sets are empty, then the indices are not defined,
i.e., $[f(u),l(u)] = \emptyset$.

In case of non-overlapping time slots and if $w(j) = u,\ j \notin \mathcal{A},$ $u$ is not empty, i.e.,
there is at least one customer $i \in \mathcal{A}$ assigned to $u$,
the following statement holds:
$$
\Theta  (j, \mathcal{A})\subseteq
\{ f(u)-1,\ldots, l(u)  \}.
$$

\subsubsection{Neighborhoods}

Our \ANSN heuristic considers two different neighborhoods for a time window $u \in \mathcal{W}$ and a tour $\mathcal{A} \in \mathcal{S}$.

\begin{compactitem}
  \item \NeighborhoodOne includes all operations with customers inside $u$ :\\
  \mbox{ $in(u, \mathcal{A}):=\{i \in \mathcal{A} :
  i \in [f(u),l(u)] \}.$}

  \item \NeighborhoodTwo represents operations with customers outside $u$ :\\ \mbox{ $out(u,\mathcal{A}):=\mathcal{A} \setminus (in(u,\mathcal{A})
  \cup \{0, n+1 \}).$}
\end{compactitem}

\begin{figure}[!ht]
  \centering
  \adjustbox{width=\textwidth}{
  %
%
\begin{tikzpicture}[every node/.style={font=\sffamily\footnotesize , line width=1pt}]
\def \lvlA {0.5cm}
\def \lvlB {-0.5cm}

\def \windL {4cm}
\def \windR {12cm}
\def \lvlTimeLine {-0.7cm}

\tikzstyle {secondWindow}=[rectangle , draw, dashed, color=red]

\tikzstyle {thirdWindow}=[rectangle, draw, dotted, color=blue]
\tikzstyle {fourthWindow}=[rectangle, draw, dotted, color=green]
\tikzstyle {firstwindow}=[rectangle, draw,  color=lightgray]

\draw[color=lightgray, line  width=0.5mm,-] (0.9cm,-0.1cm) -- (0.7cm,-0.1cm)  -- (0.7cm,1.1cm)--(0.9cm,1.1cm);
\draw[color=lightgray, line  width=0.5mm,-] (15.1cm, -0.1cm)  -- (15.3cm, -0.1cm)  -- (15.3cm, 1.1cm)--(15.1cm, 1.1cm);
\node (text) at  (0.7,-0.2cm)  {$B_{{u}}$};
\node (text) at  (15.3,-0.2cm)  {$E_{{u}}$};
\draw[color=lightgray, line  width=0.5mm, dotted, color=blue] (5.7cm, 0cm)  -- (5.9cm, 0cm) -- (5.9cm, 1cm)-- (5.7cm,1cm);
\node (text) at  (5.9,-0.2cm)  {$E_{{u-1}}$};
\draw[color=lightgray, line  width=0.5mm,dashed, color= red] (0.95cm,0cm) -- (0.75cm,0cm)  -- (0.75cm,1cm)--(0.95cm,1cm);
\draw[color=lightgray, line  width=0.5mm,dashed, color=red] (8.8cm, 0cm)  -- (9cm, 0cm)  -- (9cm, 1cm)--(8.8cm,1cm);
\node (text) at  (0.7,-0.4cm)  {$B_{{u+1}}$};
\node (text) at  (9,-0.2cm)  {$E_{{u+1}}$};
\draw[color=lightgray, line  width=0.5mm, dotted,color=green] (10.2cm,0cm) -- (10cm,0cm) -- (10cm,1cm)--(10.2cm,1cm);
\node (text) at  (10,-0.2cm)  {$B_{{u+2}}$};

\node [	secondWindow,xshift=1.8cm,yshift= {\lvlA},	minimum width=0.8cm,		minimum height=0.7cm,] (a1) {${f({u})}$};
\node [	draw,xshift=3.2cm,yshift= {\lvlA},	minimum width=0.8cm,		minimum height=0.7cm,] (a2) {${f({u})+1}$};
\node [	thirdWindow,xshift=5cm,yshift= {\lvlA},	minimum width=0.8cm,		minimum height=0.7cm,] (a3) {${f({u})+2}$};
\node [	draw,xshift=6.5cm,yshift= {\lvlA},	minimum width=0.8cm,		minimum height=0.7cm,] (a4) {$$};
\node [	secondWindow,xshift=8.2cm,yshift= {\lvlA},	minimum width=0.8cm,		minimum height=0.7cm,] (a5) {$$};
\node [	draw,xshift=9.5cm,yshift= {\lvlA},	minimum width=0.8cm,		minimum height=0.7cm,] (a6) {$$};
\node [ fourthWindow,xshift=11.4cm,yshift= {\lvlA},	minimum width=0.8cm,		minimum height=0.7cm,] (a7) {${l(u)-1}$};
\node [	draw,xshift=12.9cm,yshift= {\lvlA},	minimum width=0.8cm,		minimum height=0.7cm,] (a8) {${l(u)}$};
\node [	fourthWindow,xshift=14.3cm,yshift= {\lvlA},	minimum width=0.8cm,		minimum height=0.7cm,] (a9) {${l(u)+1}$};

\draw[decorate, ->] (a1.east) -> (a2.west) ;
\draw[decorate, ->] (a2.east) -> (a3.west) ;
\draw[decorate, ->] (a3.east) -> (a4.west) ;
\draw[decorate, ->] (a5.east) -> (a6.west) ;
\draw[decorate, ->] (a6.east) -> (a7.west) ;
\draw[decorate,->] (a7.east) -> (a8.west) ;
\draw[decorate, ->] (a8.east) -> (a9.west) ;
\draw[decorate, ->] (a4.east) -> (a5.west) ;

\draw[decorate, decoration=snake] (0.2cm,{\lvlA}) -> (a1.west) ;
\draw[decorate, decoration=snake] (15.8cm,{\lvlA}) -> (a9.east) ;


%
\node [	thirdWindow, anchor=west] at  (0cm,2cm) (t1) {$u-1$: 08:00-09:30};

\node [	secondWindow ,anchor=west] at  (3cm,2cm) (t2) { $u+1$: 09:00-09:45};

\node [firstwindow, anchor=west] at  (6cm,2cm) (t3) { $u$: 09:00-11:00};

\node [	fourthWindow, anchor=west] at  (8.5cm,2cm) (t4) {$u+2$: 10:00-11:30};

\draw[draw,-] ([shift={(-0.07,-0.7)}]a1.south  west) --
([shift={(-0.07,-0.3)}]a1.south  west);
\node [anchor=west] at ([shift={(-0.1,-0.55)}] a1.south west)  {
 {\tiny
$\Theta^{-}(j,\mathcal{A})$}};

\draw[draw, -] ([shift={(0.07,-0.7)}]a9.south  east) --
([shift={(0.07,-0.3)}]a9.south  east);
\node [anchor=west] at ([shift={(0.0,-0.55)}] a9.south east)  {
	{\tiny
		$\Theta^{+}(j,\mathcal{A})$}};


\draw[color= darkgray, draw,line  width=0.7mm, ->,>=stealth] (0cm,{\lvlTimeLine}) -- (16cm, {\lvlTimeLine});
\draw (8cm,{\lvlTimeLine -0.2cm}) node {Time};
\end{tikzpicture}}
  \caption{Illustration of the first $f(u)$ and
  last customer $l(u)$ of  $u$ in case of overlapping time windows.
  Moreover, we indicate the positions of the insertion points $\Theta^{-}(j, \mathcal{A})$
  and $\Theta^{+}(j, \mathcal{A})$ (vertical lines).
  We  notice that $\Theta  (j, \mathcal{A})$ and $
  \{ f(u)-1,\ldots, l(u)  \}$ differ as also the insertion of
  $j$ after customer
  ${l(u)+1}$ must be considered.
  Further, we observe that customer ${f(u)+2}$ is
  inside $u$ although being assigned to time window $u-1$.
  }
  \label{fig:insPoints}
\end{figure}
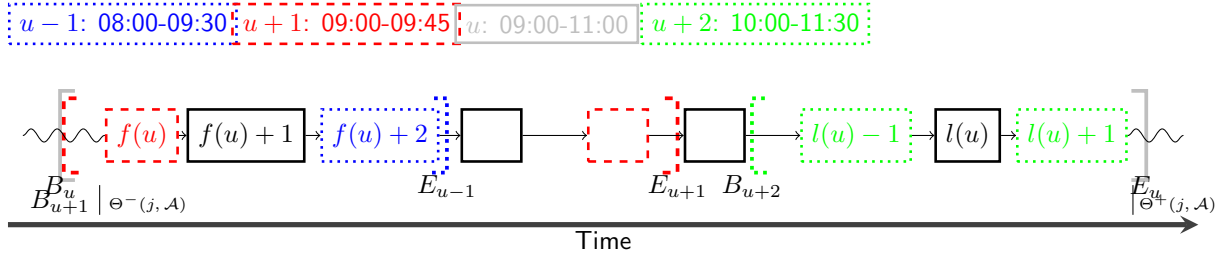

\noindent The inside of $u$ consists of customers $i \in \mathcal{A}$  that are
\begin{compactitem}
  \item assigned to time window ${u} = w(i)$,
  \item  assigned to a time window that is included in  $u$ :
  \mbox{$s_{u} \leq s_{w(i)} \leq e_{w(i)} \leq e_{u}$,}
  \item or, captured by customers of $u$ (or its included time windows),
  e.g., there exist two customers
  $j,k \in \mathcal{A}$ with $w(j)=w(k)=u$
  such that $j<i<k$ and  $i,j,k\in[n]$.
\end{compactitem}
Clearly, $in(u,\mathcal{A})$ is dependent on the actual tour sequence.
However,
in case of non-overlapping time windows,
the customers inside $u$ are exactly those who are
assigned to $u$, i.e.,
$in(u, \mathcal{A})=\{i\in \mathcal{A} : w(i)=u\}$.
In Figure \ref{fig:insPoints}, we illustrate the definitions that have been
introduced so far.

\subsubsection{Loss and Free Time}
For a tour $\mathcal{A}$ and a new customer $j$
who has to be inserted into time slot $u$, we define
\begin{align*}
  \chi_{u}^{-}(j,\mathcal{A})&:= \max
  \bigg(
  e_{f(u)}, ~
  \max_{i\in \Theta  (j, \mathcal{A}), ~ i \leq f(u)}
  \tilde{e}_{j, i} \bigg) -B_{u} ,
  \\
  \chi_{u}^{+}(j,\mathcal{A})&:= E_{u} -
  \min \bigg(
  \ell_{l(u)}, ~
  \min_{i\in \Theta  (j,\mathcal{A}), ~ i \geq l({u})}
  \tilde{\ell}_{j, i}
  \bigg).
\end{align*}
The $\chi_{u}^{-}$ corresponds to the amount of time that is
``lost'' at the beginning of time slot ${u}$.
This can be caused by
the service time required for the last customer order before (outside) $u$
or the travel time needed for going from that customer
to the first customer inside $u$.
Similarly, $\chi_{u}^{+}$
corresponds to the loss of time at the end of  time slot ${u}$
caused by the time required for traveling to the
first customer after (outside) $u$
or the service time at the last customer inside $u$.

Further, we denote $\chi_{u} (j,\mathcal{A}):=\chi_{u}^{-}
(\mathcal{A}, j) + \chi_{u}^{+} (j,\mathcal{A})$ as
the \textit{loss time} of time window ${u}$.
In case that $[f(u), l(u)]= \emptyset$,  the \textit{loss time} is given by
$$\chi_{u}(j,\mathcal{A})=
\max_{i\in \Theta  (j, \mathcal{A})}
\big(
(\tilde{e}_{j, i} -s_{u}
)
+
(
e_{u}
-
\tilde{\ell}_{j, i}
)
\big)
.
$$

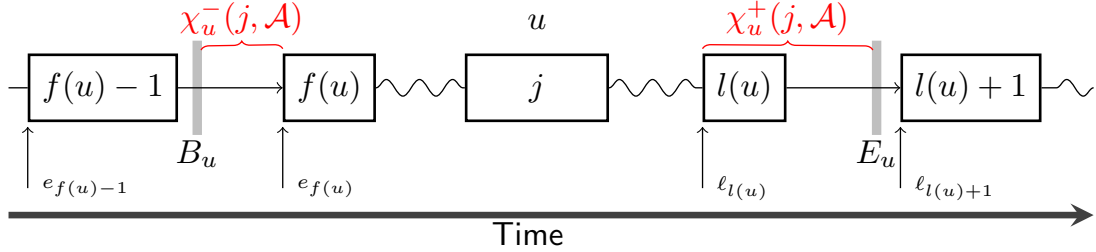
\begin{figure}[!ht]
  \centering
  \adjustbox{width=0.9\textwidth}{
  \begin{tikzpicture}[every node/.style={font=\sffamily\footnotesize, line width=0.8pt}]
\def \lvlA {0.5cm}
\def \lvlB {-0.5cm}

\def \windL {2cm}
\def \windR {9.2cm}
\def  \shiftdwn{-0.7cm}
\def \lvlTimeLine {-0.85cm}
\def \lvlBorders{-0.2cm}

\draw[color=lightgray, line  width=1mm,-] ({\windL},0cm)  -- ({\windL},1cm);

\draw[color=lightgray, line  width=1mm,-] ({\windR},0cm)  -- ({\windR},1cm);

\node [	draw,xshift=1cm,yshift= {\lvlA},	minimum width=0.8cm,		minimum height=0.7cm,] (a1) {${f(u)-1}$};
\node [	draw,xshift=3.4cm,yshift= {\lvlA},	minimum width=0.8cm,		minimum height=0.7cm,] (a2) {${f({u})}$};
\node [	draw,xshift=5.6cm,yshift= {\lvlA},	minimum width=1.5cm,		minimum height=0.7cm,] (a3) {$j$};
\node [	draw,xshift=7.8cm,yshift= {\lvlA},	minimum width=0.8cm,		minimum height=0.7cm,] (a4) {${l({u})}$};
\node [	draw,xshift=10.2cm,yshift= {\lvlA},	minimum width=0.8cm,		minimum height=0.7cm,] (a5) {${l({u})+1}$};

\draw[draw, ->] (a1.east) -> (a2.west) ;
\draw[draw, ->] (a4.east) -> (a5.west) ;
\draw[decorate,decoration=snake] (a2.east) -> (a3.west) ;
\draw[decorate,decoration=snake] (a3.east) -> (a4.west) ;

\draw[decorate,decoration=snake] (0cm,{\lvlA}) -> (a1.west) ;
\draw[decorate,decoration=snake] (11.5cm,{\lvlA}) -> (a5.east) ;

\node (bord) at (7.94cm,{\lvlA }) {};

\node (window) at (7.94cm,{\lvlA }) {};

\draw[decorate,decoration=brace, color=red] ({\windL +0.04cm}, {\lvlA+0.4cm}) -- ([shift={(0,0.05)}] a2.north  west)  node[above, midway ] {$	\chi_{u}^{-}(j,\mathcal{A})$} (0,0);

\draw[decorate,decoration=brace, color=red] ([shift={(0,0.05)}] a4.north  west) -- ({\windR - 0.04cm}, {\lvlA+0.4cm})  node[above, midway ] {$	\chi_{u}^{+}(j,\mathcal{A})$} (0,0);

\node (text) at  ({\windL},{\lvlBorders})  {$B_{{u}}$};
\node (text) at  ({\windR},{\lvlBorders})  {$E_{{u}}$};

\draw[draw, ->] ([shift={(0,-0.7)}]a1.south  west) --
([shift={(0,-0.05)}]a1.south  west);
\node [anchor=west] at ([shift={(0.0,\shiftdwn)}] a1.south west)  {
 {\tiny
$e_{f({u})-1}$}};

\draw[draw, ->] ([shift={(0,-0.7)}]a2.south  west) --
([shift={(0,-0.05)}]a2.south  west);
\node [anchor=west] at ([shift={(0.0, \shiftdwn)}] a2.south west)  {
 {\tiny
$e_{f({u})}$}};

\draw[draw, ->] ([shift={(0,-0.7)}]a4.south  west) --
([shift={(0,-0.05)}]a4.south  west);
\node [anchor=west] at ([shift={(0.0,\shiftdwn)}] a4.south west)  {
 {\tiny
$\ell_{l({u})}$}};

\draw[draw, ->] ([shift={(0,-0.7)}]a5.south  west) --
([shift={(0,-0.05)}]a5.south  west);
\node [anchor=west] at ([shift={(0.0, \shiftdwn)}] a5.south west)  {
 {\tiny
$\ell_{l({u})+1}$}};

\draw  ({0.5*\windL + 0.5*\windR},+1.2cm) node {${u}$};

\draw[color= darkgray, draw,line  width=0.7mm, ->,>=stealth] (0cm,{\lvlTimeLine}) -- (11.5cm, {\lvlTimeLine});
\draw (5.5cm,{\lvlTimeLine -0.2cm}) node {Time};
\end{tikzpicture}}
  \caption{Illustration of the loss time
  $\chi_u(j,\mathcal{A})$.
  Here, the case that the new customer
  $j$ is inserted into a tour, with non-overlapping time windows,
  between two other customers assigned to time slot ${u}$, is considered.}
  \label{fig:lossTime}
\end{figure}

Figure~\ref{fig:lossTime} illustrates
the	$\chi_{u}^{-}(j,\mathcal{A})$  and $	\chi_{u}^{+}(j,\mathcal{A})$ on a tour with non-overlapping time windows. Clearly, if $\chi_{u} (j,\mathcal{A})=0$,
then a violation of \tfeas for
$j$ can only be repaired by removing (exchanging) customers
that are inside $u$.

\begin{figure}[!ht]
  \centering
  \adjustbox{width=0.9\textwidth}{
  \begin{tikzpicture}[every node/.style={font=\sffamily\footnotesize, line width=0.8pt}]
\def \levelTime { -0.5cm}
\def \lvlA {0.5cm}
\def \windL {0.5cm}
\def \windR {10cm}

\draw[color=lightgray, line  width=1mm,-] (\windL,1cm)  -- (\windL,0cm);
\draw[color=lightgray, line  width=1mm,-] (\windR,1cm)  -- (\windR,0cm);
\node [	draw,		xshift=3.5cm, yshift= {\lvlA},	minimum width=0.9cm,	minimum height=0.6cm,]  (n1) {${f({u})}$};
\node [	draw,		xshift=5.5cm, yshift= {\lvlA},	minimum width=0.9cm,	minimum height=0.6cm,]  (n2) {${f({u})+1}$};
\node [	draw, xshift=8.5cm, yshift= {\lvlA},	minimum width=0.7cm,	minimum height=0.6cm,]  (n7) {${l({u})-1}$};
\node [	draw, xshift=10.4cm,	yshift= {\lvlA}, minimum width=0.9cm,	minimum height=0.6cm,]  (n8) {${l({u})}$};

\draw[draw, ->] (n1.east) -> (n2.west) ;
\draw[decorate,decoration=snake] (n2.east) -> (n7.west) ;
\draw[draw, ->] (n7.east) -> (n8.west) ;
\draw[double,thin, draw, <->, color=red] ({\windL+0.05cm},\lvlA) -- (n1) node[midway, below, color=red] {$\lambda_{u}(\mathcal{A})$};

\draw  ({0.5*\windL + 0.5*\windR},+1.2cm) node {${u}$};

\node (text) at  ({\windL},-0.2cm)  {$B_{{u}}$};
\node (text) at  ({\windR},-0.2cm)  {$E_{{u}}$};

\draw[color= darkgray, draw,line  width=0.7mm, ->,>=stealth] (0cm,\levelTime) -- (11.5cm, \levelTime);
\draw (5.5cm, {\levelTime-0.2cm}) node {Time};
\end{tikzpicture}}
  \caption{Illustration of free time $\lambda_{u} (\mathcal{A})$ for
  a single time window.} \label{fig:freeTimeDef}
\end{figure}
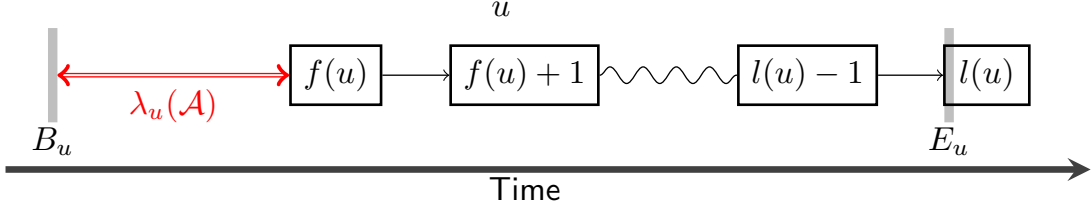

Furthermore, we want to quantify the amount of service
and travel time that is needed for inserting $j$
during time window ${u}$. For a given tour $\mathcal{A}$ the \textit{free time} of time
slot ${u}$ is defined as
\begin{align*}
  \lambda_{{u}}(\mathcal{A}):=
  \underbrace{\left( E_{u} - B_{u}  \right)}_{(\RNum{1})}
  - \underbrace{\sum\limits_{i=f({u})}^{l({u})-1} \big(s(i)  + t(i,i+1)  \big)}
  _{(\RNum{2})},
\end{align*}
where $(\RNum{1})$ is the length of ${u}$ and
$(\RNum{2})$ is the amount of service and travel time that must be
handled within ${u}$.
In case that the indices $f({u})$ and $l({u})$ are not defined,
i.e.,
$in(u, \mathcal{A})=\emptyset$,
term $(\RNum{2})$ is set to $0$.
In Figure \ref{fig:freeTimeDef}, we provide an illustration of the free time.

Considering non-overlapping time windows, the insertion of customer
$j$ after a customer
$i \in \{{f({u})-1}, \allowbreak \ldots, l({u})	\}$ requires an additional amount
of (travel and service) time that must be handled within time window ${u}$ and can be calculated by
$\lambda_{u}(\mathcal{A}) - \lambda_{{u}}(\mathcal{A} +_{i} j)$.
This assumes that all customers between
${f({u})}$ and ${l({u})}$ can be  moved arbitrarily
(while maintaining their sequence)
within time slot ${u}$, i.e.,
all customers assigned to $u$ can be seen as one consecutive block
which length is given by $(\RNum{2})$.
Hence, the above statement is weakened in the case
of overlapping time windows as some customers inside $u$ may be restricted
by their assigned time windows such that no consecutive block of customers can be formed.
Therefore, a larger amount of time than  $(\RNum{2})$ may be required
for tour $\mathcal{A}$ after the insertion of customer $j$.

\subsubsection{Feasibility and Infeasibility Conditions}

The insertion of $j$ into $\mathcal{A}$ is
infeasible if the following holds:
\begin{align}\label{cond:infeas}
  \max\limits_{i \in \Theta  (j, \mathcal{A}) } ~
  \lambda_{{u}}(\mathcal{A}  +_i j) < 0.
\end{align}

We note that Condition \eqref{cond:infeas}
solely depends on the customers inside $u$.
Hence, in case of non-overlapping time windows,
it is only dependent on the customers assigned to time slot $u$.

Moreover, in case of non-overlapping time slots,
the insertion of $j$ into $\mathcal{A}$ is feasible for at
least one insertion position if the following inequality holds:
\begin{align}\label{cond:feas}
  \max\limits_{i \in \Theta  (j, \mathcal{A}) }
  ~
  \lambda_{{u}}(\mathcal{A}  +_i j)
  -\chi_{u} (\mathcal{A}, j)  \geq 0.
\end{align}
Note that the statement does not hold in the case of overlapping time windows.

\subsubsection{Algorithm}
Finally, we can concisely describe the details of our \ANSN for the \SOPN as follows.
We temporarily assign the new customer $j$ to each time window $u \in \mathcal{W}.$
For each tour $\mathcal{A} \in \mathcal{S}_u$ we apply the steps described below,
which try to insert $j$ into $\mathcal{A}$ within its assigned time window $w(j) = u$.
If this is possible, $u$ is added to the set of available time windows  $\mathcal{T}_{j} \subseteq \mathcal{W}$.

\begin{compactitem}
  \item \mbox{\textbf{Step 1}}
  ensures that the current tour $\mathcal{A}$ fulfills \meth{CFEAS($j,\mathcal{A}$)}. .
  If $ \sum\limits_{i \in [n]} c(i) + c(j)  >  C_{\mathcal{A}}$ holds,
  then customer $a$ cannot be feasibly inserted into tour $\mathcal{A}$.
  In this case, we have to reduce the overall weight $\sum\limits_{i \in [n]} c(i)$ of $\mathcal{A}$.
  This is achieved by applying  local search operations  \onemove ~and \oneswap,
  while as few operations as possible are used.
  Therefore,  as much weight as possible is moved in each step and
  the operations stop once there is sufficient spare capacity on $\mathcal{A}$ to insert $j$, i.e.,
  $ \sum\limits_{i \in [n]} c(i) + c(j)  \leq C_{\mathcal{A}}$ holds.
  If we do not succeed in modifying $\mathcal{A}$
  such that  \meth{CFEAS($j,\mathcal{A}$)} holds, we terminate.

  \item In \mbox{\textbf{Step 2}} we aim to increase the free time $\lambda_{u}(\mathcal{A})$
  through local search operations within \NeighborhoodOne
  until the Infeasibility Condition \eqref{cond:infeas} does not hold anymore.
  $\lambda_{u}(\mathcal{A})$ is increased by applying as few operations as possible.
  That way, the  previously optimized schedule $\mathcal{S}_u$ is not  altered more than necessary.
  If a local optimum is reached,
  meaning no further improvements can be achieved by local search operations,
  and \eqref{cond:infeas} is still satisfied,
  the algorithm stops because the following steps cannot result in a feasible insertion of the new customer.

  \item \mbox{\textbf{Step 3}} is concerned with reducing the loss time $\chi_{u} (j, \mathcal{A})$ of time slot $u$ through local
  improvement operations within \NeighborhoodTwo.
  The operations are applied until either the new customer order can be inserted into the tour,
  the loss time is equal to zero $\chi_{u} (j, \mathcal{A})=0$,
  or a local optimum is reached.

  \item Finally, in \mbox{\textbf{Step 4}} we try to further increase the free time through local search operations within \NeighborhoodOne.
  The free time is increased until either the insertion of customer $j$ is possible
  or a local optimum (of the free time objective) is reached and hence,
  we are not able to insert $j$ within time slot $u$.
\end{compactitem}

Note that we apply local search operations during \mbox{\textbf{Steps 2-4}} only if \meth{CFEAS($j,\mathcal{A}$)} still holds.
Further, we apply local search operations during \mbox{\textbf{Step 3}} only if the Infeasibility Condition  \eqref{cond:infeas} does not hold for the resulting tour $\mathcal{A}$.

\subsection{Exact Approach for Solving a Subproblem} \label{ssec:exact}

In this subsection, we consider the \TSPTWN, a sub-problem of the \VRPTWN, which was first introduced by Savelsbergh \cite{savelsbergh_vehicle_1992}.
The \TSPTWN is concerned with minimizing the travel time
of a single tour $\mathcal{A} \in \mathcal{S}$ of the \VRPTWN,
while all other tours in the schedule are fixed.
Hungerl\"ander and Truden
\cite{Hungerlaender201801}
give two competitive \MILPN formulations for the \TSPTWN that
we utilize in our hybrid approaches (described in Subsection \ref{subSec: solutionApproaches}):
\begin{compactenum}
  \item A general model that can be applied to any \TSPTWN instance
  (having asymmetric travel times)
  regardless of the structure of the defined time windows $\mathcal{W}$.

  \item A second model that is tailored to the TSP with \textit{structured} Time Windows \TSPsTW.
  It is assumed that the time windows $\mathcal{W}$ are pair-wise non-overlapping and
  that the number of customers $\lvert\mathcal{C}\rvert = p$ is much larger than the number of time windows $\lvert\mathcal{W}\rvert = q$, i.e.,
  $p \gg q$, and therefore typically several customers are assigned to the same time window.
  This 	allows a simplified \MILPN formulation that performs
  significantly better.
\end{compactenum}
We refer the reader to the proceedings paper \cite{Hungerlaender201801}
for details on both \MILPs as well as a short computational study that compares
both formulations.

In contrast to the \VRPTWN the \TSPTWN is concerned with single tours.
Hence, it is unnecessary to include a capacity constraint in the \MILPN models,
as the sum of order weights of a tour is independent of the
actual sequence of the customers on the tour.

\subsection{Solution Approaches}\label{subSec: solutionApproaches}

Now we can combine
the heuristics  and \MILPs presented
above to conduct sufficiently fast
\textit{Get TWs}, \textit{Set TW}, and  \textit{Improvement steps}
in order to deal with the computational challenges that arise during the
\textit{Ordering Phase}.

\subsubsection{\textit{Get TWs step}}
In this step, we aim to quickly identify all time windows $\mathcal{T}_{j} \subseteq \mathcal{W}$ during
which a new customer $j$ can be inserted into (at least one
of) the current tours. We compare to following approaches.
\begin{compactitem}
  \item \SimpleInsertion heuristic,
  \item \ANSN heuristic,
  \item A feasibility version of the suggested \MILPs for solving the \TSPTWN, which are applied for each time slot $u \in \mathcal{W}$  	and for each $\mathcal{A} \in \mathcal{T}_u$.
\end{compactitem}

\subsubsection{\textit{Set TW step}}
Once customer $j$ has selected a time window $u \in \mathcal{T}_{j} $,
we double-check its availability in the same way as in the \textit{Get TWs step}
and then immediately insert $j$ into $u$ at the best insertion point found.
In contrast to the \textit{Get TWs step},
we run the \SimpleInsertion or the \ANSN
over all tours and all insertion points and select the feasible insertion point
that results in the least increase of the total travel time of the schedule.
Thus, we do not stop the algorithm after the first feasible insertion point is found.
In the case of the \TSPTWN  ~\MILPs we apply their standard formulations for finding an optimal tour after the insertion of $j$,
rather than their corresponding feasibility versions.

\subsubsection{\textit{Improvement step}}

In this step we aim to reduce the total travel time of the schedule.
We compare the following approaches.

\begin{compactitem}
  \item
  \onemove ~+ (\oneswap) \texttt{improvement}.
  The computationally cheap yet quite effective, \LocalSearch heuristic builds the foundation of the \textit{Improvement step}.
  We apply \onemove~  (and  \oneswap) ~ operations,
  where we focus on the \onemove ~ operations if possible,
  because they are computationally cheaper and, in general, more effective than \oneswap~ operations.
  We stop our \LocalSearch heuristic once we reach a local minimum of the
  objective function with respect to the selected neighborhood.

  \item
  \onemove ~+ \oneswap ~+ \texttt{TSP(s)TW improvement}.
  After applying the \LocalSearch heuristic we additionally run our \texttt{TSP(s)TW} \MILPN
  on	all tours that have changed since the last \textit{Improvement step}.
  We use the current tours of our delivery schedule as the initial solution for \texttt{TSP(s)TW}  \MILPN.
  While 	\onemove and \oneswap exchange customers between different tours, the
  \texttt{TSP(s)TW} \MILPN re-orders them within the tours, which makes it
  a useful complement to the \LocalSearch heuristics.
  In practice, optimizing the single  tours of a schedule to optimality has proven
  to be critical to ensure driver satisfaction as it guarantees
  that drivers do not encounter any inefficiencies on their routes.
\end{compactitem}
The above improvement procedures are arranged in ascending order with respect
to their computational effort.

During the \textit{Ordering Phase} the proposed  \LocalSearch procedures only perform improving operations.
However, the algorithms can be simply altered to a \textit{Simulated Annealing}
approach \cite{Kirkpatrick1983}
by allowing also \textit{non-improving} operations.
However, this is more suitable for the  \textit{Preparation Phase},
when more time is available for improving the delivery schedule.

\section{Computational Study and Analysis}
\label{sec:compSetup}

We want to provide a performance evaluation of the different steps and approaches
of the \textit{Ordering Phase}.
First, in Section \ref{sec:design}, we describe the design of our test instances.
In Section \ref{ssec:gettws}, we analyze how well, in terms of run time and solution quality, do the different \textit{Get TWs} approaches perform.
Then, we compare the different \textit{Improvement} approaches in
Section \ref{ssec:impr}.
Finally, in Section \ref{ssec:interplay},  we compare the performance of
different combinations of approaches for the \textit{Get TWs} and  the	\textit{Improvement step}.

\subsection{Design of the Instances}
\label{sec:design}
The  benchmark instances are designed to reflect instances
as they arise in the online grocery shopping service
of an international grocery retailer regarding
travel times, length of time windows, duration of service
times, customer order weights, and their proportions to vehicle
capacities.
Our benchmark instances are derived from those originally proposed by \cite{Cwioro2018} and have the following  characteristics.
\begin{compactitem}
  \item
  \textbf{Grid Size.}
  We consider a $20$\,km $\times$ $20$\,km square grid.
  This  roughly corresponds to  the size of a European capital such as  Vienna, Austria.
  Note that each instance corresponds to one delivery region that is served
  by one depot, which has its assigned fleet of vehicles.

  \item \textbf{Placing of Customers.}
  When placing the customers on the chosen grid, we tried to recreate
  urban settlement structures.
  Typically, these are characterized by  varying customer densities.
  To simulate this behavior, 	only \SI{20}{\percent} of the customer locations have been sampled from a two-dimensional
  uniform distribution, while the remaining \SI{80}{\percent} of the customer locations have been
  randomly assigned to $15$ clusters.
  The center (location) of each cluster $\mu =(\mu_x, \mu_y)$ is
  sampled from a two-dimensional uniform distribution.
  The shape of each cluster is defined by the covariance matrix $\Sigma=
  \begin{pmatrix}
    {\sigma_x}^2  & 0  \\
    0             &  {\sigma_y}^2
  \end{pmatrix}
  $,
  where  ${\sigma_x}^2$ and ${\sigma_y}^2$  both follow a uniform distribution.
  Furthermore, the  clusters have been rotated by a random angle between $0$ and $2\pi$.
  Customer locations have been sampled from the
  multivariate normal distribution $N(\mu, \Sigma)$ of the assigned cluster
  and
  all coordinates have been rounded to integers.
  Moreover, the numbering of the customers is randomly permuted.

  \item \textbf{Depot location.}
  We consider two different placements of the depot:
  At the center of the grid and at  the center of the top left quadrant.
  In each test setup, there are equally many instances for both variants.
  \item \textbf{Travel times.}
  We assume a travel speed of $20$\,km/h, see   Pan et al. \cite{pan2017}.
  This number can be further supported by a recent report by Vienna Public Transport
  \cite{WienerLinienEN2018}, where an average
  travel speed for their fleet of buses of $17.7$\,km/h  during the day, $17.2$\,km/h at peak times,
  and $20.0$\,km/h during evening
  hours has been reported.
  The distance between two locations is
  calculated as the Euclidean distance between them.
  Travel times are calculated proportionally to the Euclidean distances using the
  assumed travel speed and are rounded to integer seconds.
  As the Euclidean distance between two points is in general shorter than the
  shortest path distance in an underlying road network,
  we multiply all distances by a correction factor of $1.5$ \cite{cleophas_when_2014}.
  \item \textbf{Order  weights.} The order weights of customers have been
  sampled from a truncated normal distribution with
  mean of 7 and standard deviation of 2, where the lower bound is 1 and
  the upper bound is 15. The values are rounded to integers.
  Each vehicle has a loading capacity of $200$ units.
  \item \textbf{Service times.}
  We assume the service time for each customer to be 5 minutes.
  \item \textbf{Shift patterns.} All tours have the same start and end
  times. The vehicle operation times are chosen such that they
  do not restrict the actual problem.
  \item \textbf{Customer choice model.}
  Customers are iteratively inserted into the schedule
  following a simple customer choice model  that simulates the decisions
  that are usually taken by the
  customers.
  We have chosen a  simple model
  (following  \cite{cleophas_when_2014}), where every customer has just
  one desired delivery time window
  that has been set beforehand in the benchmark instance.
  If the preferred time window is not offered to the customer, we assume that
  the customer refuses to place an order.
  As in reference \cite{Agatz2011}, we assume that all delivery time windows
  are equally prominent among customers in order to obtain unbiased results that
  allow for easier identification and clearer interpretation
  of the key findings.
  Following a uniform distribution, we simulate this by a random assignment
  of each customer to one time window out of the defined
  set of time windows $\mathcal{W}$.
  \item
  For all our experiments we consider the following three sets of delivery time windows $\mathcal{W}$:
  \begin{compactitem}
    \item \textbf{Set $\windowsone$:}
    $10$ non-overlapping time windows having length of $1$ hour each,
    e.g., 08:00-09:00, 09:00-10:00,  etc.
    \item \textbf{Set $\windowstwo$:}
    $10$ overlapping time windows having length of $1.5$ hours each
    (except for  the last time window, which has $1$ hour length),
    where each window overlaps the preceding time window by $30$ minutes, e.g.,
    08:00-09:30, 09:00-10:30, etc.
    \item \textbf{Set $\windowsthree$:}	$12$ overlapping  time windows, consisting of $9$ windows having length of $1$ hour each, 08:00-09:00, 10:00-11:00,\ldots, 16:00-17:00 and (as used by Köhler et al. \cite{kohler_flexible_2019})
    $3$ time windows of $3$ hours length, morning: 08:00-11:00, noon: 11:00-14:00, afternoon: 14:00-17:00.
  \end{compactitem}

\end{compactitem}

In summary, our assumptions were chosen in order to find a good compromise between
realistic real-world instances and enabling a concise description and interpretation
of the experimental setup.

All experiments were performed on an Ubuntu 14.04 machine powered
by an Intel Xeon E5-2630V3 @ 2.4 GHz 8 core processor and 132 GB RAM.
We implemented all algorithms in Java
Version 8 and  use Gurobi 8.1.0 as \MILPN-solver in single-thread mode.
Parallelization of the applied methods is not considered.
We run each experimental configuration on $100$ instances and report average values.
Note that the absence of overlapping time windows allows using a more efficient
\MILPN~formulation of the \TSPTWN for $\windowsone$ than for  $\windowstwo$ and
$\windowsthree$ (see Section  \ref{ssec:exact}).

\subsection{Comparing 	\textit{Get TWs} Approaches}
\label{ssec:gettws}

In this section, we evaluate the performance of the proposed 	\textit{Get TWs} approaches, i.e.,
we compare the
\SimpleInsertion heuristic, the  \ANSN heuristic, and the \TSPTWInsertion approach
in terms of run time and solution quality.

\subsubsection{Test Setup}
In order to allow a proper comparison of the methods for solving the
\SOPN,
we constructed instances
which consist of
\begin{compactitem}
  \item A feasible schedule which contains $p$ customers, and
  \item a new customer order for which the availability of delivery time slots must be decided.
\end{compactitem}
To create \SOPN~instances for benchmarking we had
to create feasible  delivery schedules that are already filled with orders.
Hence, we created delivery schedules
by iteratively trying to insert $2000$ customers into each schedule.
The \SimpleInsertion heuristic was used to conduct the
feasibility checks.
The number of customers that are contained in the resulting schedule is denoted
by $\hat{p}$. We consider $\hat{p}$ being a sufficiently good approximation
of the maximal number that can be inserted into a schedule considering a given configuration.
Hence, we distinguish two scenarios.
\begin{compactenum}
  \item  In the first scenario, we perform no optimization between the insertion steps.
  \item In the second scenario, the schedule is re-optimized by applying \onemove
  after each customer insertion. This reduces the total travel time of the schedule.
  We restrict the \textit{Improvement step} to the most simple approach
  to avoid unwanted bias when evaluating the  \textit{Get TWs step}.
\end{compactenum}
In general, the schedules in the second scenario contain more orders while
utilizing the same number of vehicles.
Since the practical hardness of the \SOPN increases as
the schedules get filled up with customers, we consider \SOPN instances with different fill levels.
The \textit{fill level} $f$ of a schedule is defined as the ratio between the
number of customers $p$ in the schedule and the maximal number of
customers  $\hat{p}$ that can be inserted into the schedule.
For benchmarking at a given fill level $f$ we
select the schedule that  was generated during above described process,
containing $p=\lceil f \cdot \hat{p} \rceil$ orders.
For each generated instance we solve the \SOPN using the \SimpleInsertion
heuristic, the \TSPTWInsertion approach, and the
proposed \ANSN heuristic.
In order to investigate the differences between the considered methods,
we analyze their performance on all three  sets of time windows ($\windowsone$, $\windowstwo$, $\windowsthree$), instance sizes, and
fill levels.
Hence, we run tests on benchmark instances with $60$ tours (vehicles)
and consider fill levels of
\SI{85}{\percent}, \SI{90}{\percent}, \SI{95}{\percent}, and \SI{99}{\percent}.

\subsubsection{Results}

We report the number of feasible time slots found by each method and
required run times (mm:ss.zzz) for all scenarios, i.e., optimized and non-optimized schedules.
We report the results for  $\windowsone$ in Table \ref{tab:results_A}, for   $\windowstwo$ in Table
\ref{tab:results_C},   and
for   $\windowsthree$ in Table \ref{tab:results_D}.
Moreover, we report the number of feasible time slots
that are found by combining the findings of all three methods, denoted as \textit{Combined} in the tables.
Additionally, we report $\hat{p}$, the number of customers at \SI{100}{\percent} fill level.
All reported numbers are average values over 100 instances each.

As reference to compare against, we select the \SimpleInsertion heuristic.
Both other methods, \TSPTWInsertion and \ANSN, are entitled to
find at least the delivery time slots that are determined by the
\SimpleInsertion heuristic.
This property is guaranteed due to the construction of those methods.
Unfortunately, we can not provide an upper bound for the number of
feasible delivery time slots
as, to the best of our knowledge, there is no more powerful search method applicable for
the \SOPN in the current literature.
In our experiments, we restrict the \ANSN to \onemove operations
as preliminary experiments showed that allowing \oneswap operations yields unacceptably
long run times (up to some minutes). Also the results are only insignificantly better.
Primarily, we notice that the \SimpleInsertion heuristic returns solutions for
the \SOPN in less than one millisecond for all
considered instances.
At fill level  \SI{85}{\percent} the instances are still rather easy and hence,
already the \SimpleInsertion heuristic determines nearly all
time slots as being feasible.
While the \SimpleInsertion heuristic still performs well at \SI{90}{\percent} and
\SI{95}{\percent} on  optimized schedules,
it performs poorly on non-optimized schedules with the same fill levels
due to the lower quality of the schedules.

Further, we observe that the \TSPTWInsertion approach yields a slight improvement
over the \SimpleInsertion heuristic in terms of available time windows
at  \SI{85}{\percent} -- \SI{95}{\percent}.
However, a significant improvement can be observed when it
is applied to non-optimized schedules at  \SI{99}{\percent} fill level.
\TSPTWInsertion shows acceptable run times  for  $\windowsone$.
In contrast, run times for  $\windowstwo$ and   $\windowsthree$
are between $3$ seconds and nearly $4$ minutes and thus, unacceptable.
Hence, considering these findings, the \TSPTWInsertion
turns out to be impractical for \AHDN systems.
Moreover, we notice that the \ANSN  yields significantly more feasible
time slots than \SimpleInsertion ~(and \TSPTWInsertion) on non-optimized schedules
at  \SI{95}{\percent} and  \SI{99}{\percent}. Similar behavior can be observed on optimized schedules at  \SI{99}{\percent}.
The run times of the \ANSN stay below $1$ second for nearly all experiments with up to  \SI{95}{\percent} fill level.
The \ANSN clearly performs best in terms of solution quality on instances
having   \SI{99}{\percent} fill level resulting in up to $11$ times more available
delivery time windows than \SimpleInsertion.
However, on those instances its run time reaches up to $9$ seconds.
It is worth pointing out that the performance of the \ANSN is nearly constant over all three sets of time windows
showing that it can also deal with instances having overlapping delivery time windows.

In general, we notice a slight performance drop of the \ANSN~(compared to \TSPTWInsertion)
from optimized schedules to non-optimized schedules.
This can be explained by the fact that identifying feasible time windows is
less hard for non-optimized schedules as they contain less orders on average.
Also there is more potential for improvement when applying
\TSPTWInsertion as the single tours have not been improved in any way after inserting the customers.
Also, we observe that \textit{Combined} shows only a marginal improvement over the
\ANSN for optimized schedules.
On the other hand, we notice a strong improvement of \textit{Combined} compared to
\TSPTWInsertion and  \ANSN for non-optimized schedules at \SI{99}{\percent} fill level.
This can be explained by the fact that \TSPTWInsertion has a larger potential of finding a feasible insertion (compared to \SimpleInsertion)
if the tours are of low quality (as there is more room to rearrange the tours)
and the already observed performance drop of the \ANSN.

Our experiments show that the \ANSN heuristic is capable of finding a larger number of
feasible delivery slots than the \SimpleInsertion heuristic, requiring run times that
are suited for \AHDN services when dealing with moderately sized problem instances.
However, to efficiently tackle very large instances, parallelization of the \ANSN
is advised.

In summary, the \ANSN heuristic is clearly the best method for solving the \SOPN
when being concerned with the solution quality while
the \SimpleInsertion heuristic is the method of choice in case of very tight time restrictions.

\begin{table}[!ht]
  \scriptsize
  \caption{Summary of the computational experiments concerning the
  \textit{Get TWs} approaches for $\windowsone$
  considering non-optimized and \onemove-optimized   schedules.}
  \centering

  	\resizebox{\linewidth}{!}{
\begin{tabular}{|l|rrrr|rrrr|}
\hline
\multicolumn{9}{|c|}{\bf $\windowsone$ (10 windows) - 60 vehicles}
\\
\hline
&
\multicolumn{4}{c|}{\bf non-optimized schedules}
&
\multicolumn{4}{c|}{\bf optimized schedules}
\\
Avg. $\hat{p}$

&
\multicolumn{4}{c|}{ 408.0}
&
\multicolumn{4}{c|}{ 1907.3}
\\
Fill level
& 85\%    & 90\%  & 95\% &  99\%
& 85\%    & 90\%  & 95\% &  99\%
\\
\hline\hline
{\bf Avg. run time}
&&&&&&&&\\
{(mm:ss.zzz)}
&&&&&&&&\\
\SimpleInsertion
& 00:00.001 & 00:00.001 & 00:00.001 & 00:00.001
& 00:00.001 & 00:00.001 & 00:00.001 & 00:00.001\\

\TSPTWInsertion
& 00:00.223 & 00:00.422 & 00:01.112 & 00:01.970
& 00:00.592 & 00:00.580 & 00:00.743 & 00:00.455 \\

\ANSN
& 00:00.024 & 00:00.089 & 00:00.468 & 00:04.031
& 00:00.049 & 00:00.052 & 00:00.077 & 00:01.377 \\

\hline
{\bf Avg. number}
&&&&&&&&\\
{\bf of feasible slots}
&&&&&&&&\\
\SimpleInsertion
& 9.82 & 9.47 & 7.67 & 2.00
& 9.98 & 9.98 & 9.83 & 4.20
\\
\TSPTWInsertion
& 9.85 & 9.62 & 8.60 & 4.30
& 9.98 & 9.98 & 9.83 & 4.29
\\
\ANSN
& 9.96 & 9.94 & 9.93 & 9.39
& 9.98 & 9.98 & 9.98 & 9.63
\\
\hline
Combined
& 9.96 & 9.94 & 9.93 & 9.43
& 9.98 & 9.98 & 9.98 & 9.63
\\ \hline

\end{tabular}
}

  \label{tab:results_A}
\end{table}

\begin{table}[!ht]
  \scriptsize
  \caption{Summary of the computational experiments  concerning the
  \textit{Get TWs} approaches  for  $\windowstwo$
  considering non-optimized and \onemove-optimized  schedules.}
  \centering

  \resizebox{\linewidth}{!}{
\begin{tabular}{|l|rrrr|rrrr|}
\hline
\multicolumn{9}{|c|}{\bf $\windowstwo$ (10 windows) - 60 vehicles}
\\
\hline
&
\multicolumn{4}{c|}{\bf non-optimized schedules}
&
\multicolumn{4}{c|}{\bf optimized schedules}
\\
Avg. $\hat{p}$
&
\multicolumn{4}{c|}{ 628.0}
&
\multicolumn{4}{c|}{ 1880.0}

\\
Fill level
& 85\%    & 90\%  & 95\% &  99\%
& 85\%    & 90\%  & 95\% &  99\%
\\
\hline\hline
{\bf Avg. run time}
&&&&&&&&\\
{(mm:ss.zzz)}
&&&&&&&&\\
\SimpleInsertion

& 00:00.001 & 00:00.001 & 00:00.001 & 00:00.001
& 00:00.001 & 00:00.001 & 00:00.001 & 00:00.001
\\
\TSPTWInsertion
& 00:03.505 & 00:11.182 & 00:52.761 & 02:22.804
& 00:41.024 & 01:01.611 & 01:03.537 & 00:11.598
\\
\ANSN
& 00:00.032 & 00:00.138 & 00:00.683 & 00:07.837
& 00:00.036 & 00:00.040 & 00:00.066 & 00:08.343
\\

\hline
{\bf Avg. number}
&&&&&&&&\\
{\bf of feasible slots}
&&&&&&&&\\
\SimpleInsertion
& 9.57 & 9.05 & 7.40 & 0.76
& 9.98 & 9.98 & 9.68 & 0.72
\\
\TSPTWInsertion
& 9.84 & 9.68 & 8.96 & 5.04
& 10.00 & 10.00 & 9.84 & 0.70
\\
\ANSN
& 9.99 & 9.95 & 9.83 & 8.43
& 10.00 & 10.00 & 9.84 & 6.36
\\
\hline
Combined
& 9.99 & 9.95 & 9.87 & 8.76
& 10.00 & 10.00 & 10.00 & 6.36
\\ \hline

\end{tabular}
}

  \label{tab:results_C}
\end{table}

\begin{table}[!ht]
  \scriptsize
  \caption{Summary of the computational experiments  concerning the
  \textit{Get TWs} approaches for  $\windowsthree$ considering non-optimized  and \onemove-optimized schedules.}
  \centering
  	\resizebox{\linewidth}{!}{
\begin{tabular}{|l|rrrr|rrrr|}
\hline
\multicolumn{9}{|c|}{\bf  $\windowsthree$ (12 windows)  - 60 vehicles}
\\
\hline
&
\multicolumn{4}{c|}{\bf non-optimized schedules}
&
\multicolumn{4}{c|}{\bf optimized schedules}
\\
Avg. $\hat{p}$
&
\multicolumn{4}{c|}{ 406.0}
&
\multicolumn{4}{c|}{ 1897.9}
\\
Fill level
& 85\%    & 90\%  & 95\% &  99\%
& 85\%    & 90\%  & 95\% &  99\%
\\
\hline\hline
{\bf Avg. run time}
&&&&&&&&\\
{(mm:ss.zzz)}
&&&&&&&&\\
\SimpleInsertion
& 00:00.001 & 00:00.001 & 00:00.001 & 00:00.001
& 00:00.001 & 00:00.001 & 00:00.001 & 00:00.001
\\
\TSPTWInsertion
& 00:02.645 & 00:15.541 & 01:27.422 & 05:20.061
& 00:56.969 & 00:54.405 & 01:04.239 & 00:17.879
\\
\ANSN
& 00:00.012 & 00:00.112 & 00:00.847 & 00:08.412
& 00:00.052 & 00:00.057 & 00:00.073 & 00:08.981
\\
\hline
{\bf Avg. number}
&&&&&&&&\\
{\bf of feasible slots}
&&&&&&&&\\
\SimpleInsertion
& 11.84 & 11.02 & 8.54 & 1.21
& 12.00 & 12.00 & 11.90 & 0.73
 \\
\TSPTWInsertion
& 11.93 & 11.73 & 10.78 & 6.65
& 12.00 & 12.00 & 11.92 & 0.69
\\
\ANSN
& 12.00 & 11.98 & 11.94 & 10.16
& 12.00 & 12.00 & 12.00 & 7.59
\\
\hline
Combined
& 12.00 & 11.99 & 11.97 & 10.73
& 12.00 & 12.00 & 12.00 & 7.59
 \\
 \hline

\end{tabular}
}

  \label{tab:results_D}
\end{table}

\subsection{Comparing \textit{Improvement} Approaches}

\label{ssec:impr}
To compare the proposed \textit{Improvement}approaches.
We perform experiments where we iteratively insert new customers into the schedule,
simulating customers placing orders online.

Due to the iterative benchmark setup, we can omit the \textit{Set TW step} and insert the
new order without double-checking the availability of the selected delivery
time slot. Again, to avoid bias, we stick to the most simple \textit{Get TWs} approach,
the \SimpleInsertion heuristic.
Then, for the \textit{Improvement step} we compute the following metrics:
\begin{compactitem}
  \item
  Average improvement
  over \textit{Insertion step}:
  The average reduction of the objective function when applying the optimization approaches
  to the schedule after inserting the new customer
  (given in percentage).
  \item
  Average improvement of the cost of \textit{Insertion}:
  The
  average reduction of the objective function relative to the
  increase of the objective function caused by inserting  the new customer
  (given in percentage).
  \item Average number of \MILPs solved for the \TSPsTWN.
  \item Average runtime of each improvement strategy.
\end{compactitem}
Additionally, we report the average total number of customer orders that have
been inserted into the final schedules.
Note that for the \MILPs we set a time limit of 60 seconds.

\subsubsection{Average-Sized Grocery Home Delivery Problems}

First, we want to analyze the \textit{Improvement} approaches
for instance sizes which we found to appear most commonly in practice.
Hence, we consider $500$ customers that are served by
$16/18/20$ vehicles having capacity of $200$ units each.
The number of used vehicles corresponds to the practical
difficulty of the instances. The numbers are chosen such that the instances are reasonably difficult.
In that sense, using $20$ vehicles results in accepting nearly all $500$ customers on average.
These instances are designed to reflect the majority of delivery regions as they
were encountered during our project with a leading supermarket chain.

The results for these experiments are reported in Table   \ref{tab:opt_average_size}.
We observe that all approaches considered are applicable in an online service as the
average run time per step is below $4$ seconds, which is very reasonable for instances of
this size. Further, a reduction of our objective function by
\SI{0.29}{\percent} to  \SI{0.60}{\percent}
per step is remarkable as between two \textit{Improvement steps} the schedule is
altered only by the insertion of one customer.
This can be further underlined by the reported average reduction of the cost of
inserting the new customer ranging from
\SI{35.82}{\percent} to  \SI{63.48 }{\percent}.
These numbers show that our approaches meet the requirements of modern \AHDN systems.
The experiments reveal that  \onemove~+~\oneswap clearly outperforms \onemove in terms of
improving the objective function (across all three types of time windows). Also solving the \TSPTWN afterwards results in further
improvement of the objective.
Considering different delivery time windows, we notice that the approaches perform
best with respect to runtime on instances with $\windowsone$ and worst on instances
with $\windowstwo$.
While there is a slight difference for the \LocalSearch operations,
the difference is nearly 3 seconds when additionally applying the \TSPTWN~\MILPN.
This is due to the fact that the
absence of overlapping time windows allows for a more efficient \MILPN~formulation (Section \ref{ssec:exact}).
Thus, during peak times and in case of overlapping time windows we advise to
stick to \onemove~ (or \onemove~+~\oneswap).
Further, the average number of customers that can be inserted into the schedule
deviates at most by $9$ (+\SI{1.9}{\percent}) between $\windowsone$ and $\windowstwo$.
Hence, allowing overlapping time windows accounts for a small benefit concerning
the degree of capacity utilization.
Similarly, in case of overlapping time
windows ($\windowstwo$ and $\windowsthree$) the travel time reduction is slightly
larger than for $\windowsone$.

\subsubsection{Dealing with large Problem instances - \textit{Improvement step} only every $i$th iteration}

Large supermarket chains  offer their services across the whole country.
Caused by the different geographies, the sizes of the delivery regions that
are covered by a depot strongly vary ranging from a few hundred  up to around $2000$ customers per day.
Dealing with such large delivery regions is especially
challenging  during periods of many customer request coming in within a short time frame.
To accommodate these periods of high request frequency,
we propose to run the \textit{Improvement step} only after each $i$th successful
\textit{Insertion step}, instead of after each.

We want to validate this idea by running computational experiments.
We consider instances with $2000$
customers  (the largest number we encountered in practice) and
$80$ vehicles for these experiments and report the results in Table \ref{tab:opt_average_large}.
Each column shows results for different values of $i$.
For these experiments we can only report the average
improvement of the schedule when applying the respective improvement strategy.
First, we notice that the runtime of the \textit{Improvement step} increases with $i$:
the more often we skip an \textit{Improvement step}, the longer it takes to improve the schedule's total travel time.
Moreover, we observe an increased
improvement per step with increasing $i$.
Apparently, the improvement that was omitted can be made up (up to a certain extent),
by applying the \textit{Improvement step} at a later point in time.
It shows that performing the \textit{Improvement step}
only after every $i$th successful insertion is a viable option.
While  \onemove stays below $6$ seconds for $i=10$, its runtime increases up to $13$ seconds for $i=30$.
The runtimes of  \onemove~+~\oneswap are between $1$ minute and $2$ minutes (overlapping time windows), which
is still acceptable. We observe that running the \TSPTWN~\MILPN does increase the runtimes (on top of  \LocalSearch)
only insignificantly while still showing some additional improvements of the objective.
In general, the observations from the previous experiments with $500$ customers carry over to this experiments.
Again we notice shorter runtimes for $\windowsone$ than for $\windowstwo$ and $\windowsthree$.
However, we observe less significant  differences  than for the experiments with $500$ customers.

In summary, the results show that skipping the
\textit{Improvement step} allows us to deal with temporarily high customer request rates,
even for very large schedules with a vast number of customers.
Furthermore, we see that applying the \textit{Improvement step} less often leads to an increased
improvement per step at the cost of longer runtimes.
Finally, note that  triggering the \textit{Improvement step} dynamically
when there are no new requests is also a valid option.

\begin{table}[!ht]
  \scriptsize
  \caption{
  Summary of the computational experiments for the  \textit{Improvement} approaches considering instances with $500$ customers.}
  \centering
  \resizebox{\linewidth}{!}{
\begin{tabular}{|l|rrr|rrr|rrr|}
\hline
\multicolumn{10}{|c|}{\bf Average-sized grocery home delivery problems}
\\
\hline
&
\multicolumn{3}{c|}{\bf $\windowsone$}
&
\multicolumn{3}{c|}{\bf $\windowstwo$}
&
\multicolumn{3}{c|}{\bf $\windowsthree$}
\\
Tours:
& 16 & 18 & 20
& 16 & 18 & 20
& 16 & 18 & 20\\
\hline
Avg. $\hat{p}$
&450.80  & 491.60 & 498.00
 & 459.60  & 495.80 & 499.20
 & 453.00  & 496.60 & 499.00 \\

Avg. number of &&&&&&&&&\\
time windows offered
& 9.01   &  9.84   & 9.97
& 9.21      & 9.93    & 9.98
& 10.83     & 11.92    & 12.00
\\
\hline\hline
{\bf Avg. run time}
 &&&&&&&&&\\
{(mm:ss.zzz)}
 &&&&&&&&&\\
 \texttt{1move}
 & 00:00.015   & 00:00.099    & 00:00.108
 & 00:00.097   & 00:00.124    & 00:00.134
 & 00:00.090   & 00:00.111    & 00:00.120
\\
\texttt{1move+1swap}
 & 00:00.632   & 00:00.776    & 00:00.766
 & 00:00.952   & 00:01.170    & 00:01.145
 & 00:00.925   & 00:01.106    & 00:01.106
\\
 \texttt{1move+1swap+TSP(s)TW}
 & 00:00.678   & 00:00.833    & 00:00.821
 & 00:03.481   & 00:03.757    & 00:03.655
 & 00:01.314   & 00:01.478    & 00:01.485
\\
\hline
{\bf Avg. improvement}
 &&&&&&&&&\\
{\bf over Insertion (\%)}
 &&&&&&&&&\\
\texttt{1move}
 & 0.30  & 0.29     & 0.29
 & 0.37  & 0.37     & 0.37
 & 0.34  & 0.34     & 0.34
\\
\texttt{1move+1swap}
 & 0.41  & 0.40     & 0.39
 & 0.52  & 0.50     & 0.50
 & 0.50  & 0.50     & 0.49
\\
 \texttt{1move+1swap+TSP(s)TW}
 & 0.45 & 0.43     & 0.42
 & 0.60  & 0.57     & 0.56
 & 0.55  & 0.54     & 0.53
\\
\hline
{\bf Avg. improvement}
 &&&&&&&&&\\
{\bf of cost of Insertion (\%)}
 &&&&&&&&&\\
 \texttt{1move}
 & 35.82  & 38.00     & 38.76
 & 38.36  & 40.75     & 41.66
 & 37.11  & 39.39     & 39.58
\\
 \texttt{1move+1swap}
 & 49.88  & 51.34     & 51.44
 & 53.62  & 55.37    & 55.64
 & 54.65  & 57.62     & 57.23
\\
 \texttt{1move+1swap+TSP(s)TW}
 & 53.83  & 55.34    & 55.26
 & 61.80  & 63.33     & 63.48
 & 59.53  & 62.37     & 61.90
\\
\hline
{\bf Avg. number of}
 &&&&&&&&&\\
{\bf \MILPs solved}
 &&&&&&&&&\\
 \texttt{1move+1swap+TSP(s)TW}
 & 2.01  & 2.09     & 2.11
 & 2.36  & 2.49     & 2.48
 & 2.32  & 2.42     & 2.42
\\
\hline
\end{tabular}
}

  \label{tab:opt_average_size}
\end{table}
\begin{table}[!ht]
  \scriptsize
  \caption{Summary of the computational experiments for the \textit{Improvement} approaches
  considering instances with $2000$ customers served by $80$ vehicles.
  The \textit{Improvement step} is triggered  after every $i$th ($i=10,20,30$) successful insertion of a new customer.}
  \centering
  \resizebox{\linewidth}{!}{
\begin{tabular}{|l|rrr|rrr|rrr|}
  \hline
  \multicolumn{10}{|c|}{\bf Large-scale problem instances - 80 vehicles}
  \\
  \hline
  &
  \multicolumn{3}{c|}{\bf $\windowsone$}
  &
  \multicolumn{3}{c|}{\bf $\windowstwo$}
  &
  \multicolumn{3}{c|}{\bf $\windowsthree$}
  \\
  i:
  & 10 & 20 & 30
  & 10 & 20 & 30
  & 10 & 20 & 30\\
  \hline
  Avg. $\hat{p}$
  & 1997.20  & 1998.50 & 1998.50
  & 1999.40  & 1999.40 & 1999.40
  & 1999.80  & 1999.80 & 1999.80  \\

  Avg. number of &&&&&&&&&\\
  time windows offered
  & 9.99  & 9.99    & 9.99
  & 9.99      & 9.99    & 9.99
  & 9.99     & 9.99    & 9.99
  \\
  \hline\hline
  {\bf Avg. run time}
  &&&&&&&&&\\
  {(mm:ss.zzz)}
  &&&&&&&&&\\
  \texttt{1move}
  & 00:04.303   & 00:06:907   & 00:10:578
  & 00:05.360   & 00:08.930    & 00:12.587
  & 00:04.939   & 00:08.126    & 00:11.208
  \\
  \texttt{1move+1swap}
  & 00:59.003   & 01:28:161    & 01:57.892
  & 01:25.017   & 02:02.721    & 02:31.010
  & 01:13.358   & 01:40.780    & 02:10.395
  \\
  \texttt{1move+1swap+TSP(s)TW}
  & 00:59.980   & 01:29.274    & 01:59.241
  & 01:40.952   & 02:25.107    & 02:56.644
  & 01:16.316   & 01:50.400    & 02:14.499
  \\
  \hline
  {\bf Avg. improvement}
  &&&&&&&&&\\
  {\bf over Insertion (\%)}
  &&&&&&&&&\\
  \texttt{1move}
  & 1.62  & 3.27    & 4.93
  & 1.69  & 3.54     & 5.12
  & 1.60  & 3.29     & 4.82
  \\
  \texttt{1move+1swap}
  & 2.29  & 4.16     & 5.94
  & 2.46  & 4.57     & 6.28
  & 2.33  & 4.26    & 5.92
  \\
  \texttt{1move+1swap+TSP(s)TW}
  & 2.35  & 4.24    & 6.05
  & 2.56  & 4.73     & 6.48
  & 2.39  & 4.35    & 6.03
  \\
  \hline

  \hline
  {\bf Avg. number of}
  &&&&&&&&&\\
  {\bf \MILPs solved}
  &&&&&&&&&\\
  \texttt{1move+1swap+TSP(s)TW}
  & 23.43  & 30.64     & 33.58
  & 24.65  & 31.41     & 34.48
  & 24.09  & 31.09     & 34.22
  \\
  \hline
\end{tabular}
}

  \label{tab:opt_average_large}
\end{table}

\subsection{Interplay of \textit{Get TWs} and	\textit{Improvement} Approaches}
\label{ssec:interplay}

In this final experiments, we want to find out which combinations of
the different approaches for the \textit{Get TWs}- and  the	\textit{Improvement step}
are most beneficial
and which should be avoided.

From Section \ref{ssec:gettws}   we learn that
\SimpleInsertion is the fasted method for the \textit{Get TWs step} showing a solid performance,
while  \ANSN is the best method in terms of solution quality having the draw back
that it can only be applied when the customer request rate is moderate (or instances are small).
In Section \ref{ssec:impr}  we observe that
\onemove   is  a solid approach  for the 	\textit{Improvement step}
that also scales   well  for large problem instances.
The  application of more sophisticated \LocalSearch operations in combination with
an exact approach for a selected sub-problem (\onemove ~+ \oneswap ~+ \texttt{TSP(s)TW improvement})
shows the best performance in terms of solution quality at the price of high (but still acceptable) run times.

Further evaluations are based on the average-sized grocery delivery use case (Section \ref{ssec:impr})
and will focus on aforementioned   \textit{Get TWs} and 	\textit{Improvement} approaches.
Hence, we compare the resulting four combinations
$
\{
\SimpleInsertion, ~\ANSN
\} $
$\times$ $
\{
\onemove, ~\onemove ~+ \oneswap ~+ \texttt{TSP(s)TW improvement}\}
$
concerning the following key figures:
\begin{compactitem}
  \item Average run time of the  \textit{Get TWs}- and 	\textit{Improvement step}.

  \item Average number of offered delivery time windows.

  \item Average total number of customer orders that have
  been inserted into the final schedules.
\end{compactitem}

\subsubsection{Results}

In Table \ref{tab:comb} we report the results of this experiments.
First, we notice that now when analyzing the interplay of
the get  \textit{GetTWs step} and the \textit{Improvement step}
the differences between   \ANSN  and \SimpleInsertion become less evident.
\ANSN only show  little benefit compared to \SimpleInsertion in terms of the
number accepted orders $\hat{p}$ (at most \SI{0,4}{\percent} improvement)
and the number of offered time windows.
The use of \ANSN  reduces
the run time of the \textit{Improvement step}.
This effect is most evident when overlapping time windows are used ($\windowstwo$ and $\windowsthree$),
especially for the  \onemove ~+ \oneswap+ ~\texttt{TSP(s)TW} approach where a reduction of the run time of up to \SI{71,4}{\percent} is observed.
Presumably, \ANSN creates better schedules when inserting the new customer,
and therefore the \textit{Improvement} approaches have a better starting solution.

In summary,  we can conclude that using the \ANSN for the \textit{Get TWs} (and \textit{Set TWs}) \textit{step}
is preferred as long as the instances are sufficiently small (as in our use case) such
that the  \textit{GetTWs step} can be performed in accordance with the run time requirements
of the considered application of the \AHDN service.
However, the \ANSN  gives a slight improvement compared to
\SimpleInsertion ~(as already shown in \ref{ssec:impr}) and should therefore
be utilized if the frequency of incoming orders is sufficiently low, such that the
additionally required run times do not cause issues.
If the expected time between incoming order requests temporarily increases, e.g.,
during peak times,
one can switch to the \onemove heuristic without having to fear major drawbacks.

\begin{table}[!ht]
  \caption{
  Summary of the computational experiments
  for different combinations of  \textit{Get Tws} approaches
  ($\SimpleInsertion, ~\ANSN$)
  and
  \textit{Improvement} approaches
  ($\onemove, ~\onemove ~+~\oneswap ~+~\texttt{TSP(s)TW improvement}$)
  considering instances with $500$ customers.
  }
  \centering
  \scriptsize
  \resizebox{\linewidth}{!}{

\begin{tabular}{|l|rrr|rrr|rrr|}
\hline
\multicolumn{10}{|c|}{\bf Average-sized grocery home delivery problems}
\\
\hline
&
\multicolumn{3}{c|}{\bf $\windowsone$}
&
\multicolumn{3}{c|}{\bf $\windowstwo$}
&
\multicolumn{3}{c|}{\bf $\windowsthree$}
\\
Tours:
& 16 & 18 & 20
& 16 & 18 & 20
& 16 & 18 & 20\\
\hline
\hline
\multicolumn{10}{|c|}{\texttt{Get TWs}: \SimpleInsertion, \texttt{Improvement}: \onemove}\\
\hline
Avg. $\hat{p}$
&453.35  & 491.80   & 498.65
&455.75  & 495.10   & 498.65
&451.00  & 495.50   & 499.60
\\

Avg. number of &&&&&&&&&\\
time windows offered
& 9.07   &  9.84    & 9.97
& 9.15  &  9.90    & 9.98
& 10.81     & 11.88   & 11.99
\\
\hline
{\bf Avg. run time}
 &&&&&&&&&\\
{(mm:ss.zzz)}
 &&&&&&&&&\\
 \texttt{Get TWs}
 & 0:00.001   & 0:00.001    & 0:00.001
 & 0:00.001   & 0:00.001    & 0:00.001
 & 0:00.001   & 0:00.001   & 0:00.001
\\
\texttt{Improvement}
& 0:00.076   & 0:00.103    & 0:00.114
& 0:00.103   & 0:00.134    & 0:00.138
& 0:00.096   & 0:00.123    & 0:00.127
\\
\hline
\hline
\multicolumn{10}{|c|}{\texttt{Get TWs}: \SimpleInsertion, \texttt{Improvement}: \onemove ~+ \oneswap ~+ \texttt{TSP(s)TW}}\\
\hline
Avg. $\hat{p}$
&454.25  & 494.00   & 499.05
&456.95  & 495.95   &  499.00
&452.80  & 496.00   & 499.65
\\

Avg. number of &&&&&&&&&\\
time windows offered
& 9.09   &  9.87    &  9.98
& 9.15   &  9.92    &  9.98
& 10.86     & 11.91   & 12.00
\\
\hline
{\bf Avg. run time}
 &&&&&&&&&\\
{(mm:ss.zzz)}
 &&&&&&&&&\\
 \texttt{Get TWs}
 & 0:00.001   & 0:00.001    & 0:00.001
 & 0:00.001   & 0:00.001    & 0:00.001
 & 0:00.001   & 0:00.001   & 0:00.001
\\
\texttt{Improvement}
& 0:00.758   & 0:00.852    & 0:00.957
& 0:03.139   & 0:03.234    & 0:03.133
& 0:01.345   & 0:01.672    & 0:01.598
\\
\hline
\hline

\multicolumn{10}{|c|}{ \texttt{Get TWs}: \ANSN, \texttt{Improvement}: \onemove}\\
\hline
Avg. $\hat{p}$
&455.35  & 498.50   & 499.55
&458.00 & 499.05   & 499.75
&453.10  & 499.44   & 499.80
\\

Avg. number of &&&&&&&&&\\
time windows offered
& 9.11   &   9.97    & 9.99
& 9.16   &  9.98    & 9.99
& 10.88     & 11.99   & 12.00
\\
\hline
{\bf Avg. run time}
 &&&&&&&&&\\
{(mm:ss.zzz)}
 &&&&&&&&&\\
 \texttt{Get TWs}
 & 0:00.098   & 0:00.022    & 0:00.017
 & 0:00.090   & 0:00.014    & 0:00.012
 & 0:00.120   & 0:00.020    & 0:00.017
\\
\texttt{Improvement}
& 0:00.060   & 0:00.077    & 0:00.086
& 0:00.073   & 0:00.091    & 0:00.100
& 0:00.073   & 0:00.092    & 0:00.098
\\
\hline
\hline

\multicolumn{10}{|c|}{ \texttt{Get TWs}: \ANSN, \texttt{Improvement}: \onemove ~+ \oneswap ~+ \texttt{TSP(s)TW}}\\
\hline
Avg. $\hat{p}$
&455.55 & 498.95   & 499.65
&457.85  & 499.15   & 499.75
&452.90  & 499.55   & 499.80
\\

Avg. number of &&&&&&&&&\\
time windows offered
& 9.12   &  9.98   & 9.99
& 9.16   &  9.99    & 9.99
& 10.87     & 11.99   & 12.00
\\
\hline
{\bf Avg. run time}
 &&&&&&&&&\\
{(mm:ss.zzz)}
 &&&&&&&&&\\
 \texttt{Get TWs}
 & 0:00.104   & 0:00.026    & 0:00.020
 & 0:00.088   & 0:00.014    & 0:00.014
 & 0:00.114   & 0:00.023    & 0:00.019
\\
\texttt{Improvement}
& 0:00.713   & 0:00.889    & 0:00.862
& 0:00.898   & 0:01.138   & 0:01.099
& 0:00.883   & 0:01.158    & 0:01.161
\\
\hline
\end{tabular}
}

  \label{tab:comb}
\end{table}

\section{Conclusion}
\label{sec:conclude}

In this work, we considered an \textit{Attended Home Delivery} \AHD system in the context of  an
online grocery shopping service offered by an international  retailer.
\AHDN systems are used whenever
the customers must be
present when their deliveries arrive.
Therefore, the  company and the
customer must both agree on a  time window
during which  the delivery can be  guaranteed.

We focused on the phase during which customers place their orders through a
web service.
Generally, this is the most challenging phase of an \AHDN system
from a computational point of view.
As for most \AHDN approaches in the literature, we considered
a Vehicle Routing Problem with
Time Windows \VRPTW as the underlying optimization problem.
The online characteristic of this phase requires that the delivery schedule is
built dynamically as  new orders are placed.
We split the computations  into four steps and
proposed solution approaches that allow to
determine which delivery time windows can be
offered to potential customers and to iteratively build the schedule.

Finally, we presented a comprehensive experimental evaluation of the proposed heuristic approaches,
which are based on \LocalSearch operations and \textit{Mixed-Integer Linear Programming} formulations.
Our goal was to determine the efficiency of the approaches on benchmark sets
motivated by an international supermarket chain's online grocery shopping service.
We elaborated certain aspects of the problem by varying the structure of the time windows,
the number of available vehicles, and the number of total customer requests.
In particular, we compare different approaches for inserting new customers into
the existing delivery schedule and for re-optimizing the schedule once a new customer has been added to the schedule.
The computational study shows that the suggested algorithms can solve the
considered benchmark instances sufficiently fast to comply with the very strict runtime
restrictions arising in \AHDN systems with high customer request rates.

\section*{Statements and Declarations}

\subsection*{Acknowledgments}
The authors would like to thank Mario Ruthmair for his valuable advice and feedback.

\subsection*{Funding} This work was supported by Lakeside Labs GmbH, Klagenfurt, Austria and funding from the European
Regional Development Fund and the Carinthian Economic Promotion Fund (KWF) under grant 20214/31942/45906.
Christian Truden is supported by a  grant of the Government of Carinthia within the CARINTHIja 2020 project.

\subsection*{Competing interests}
The authors have no relevant financial or non-financial interests to disclose.

\subsection*{Availability of data and materials}

The benchmark instances are available from \url{https://bit.ly/3doVvve}.

\end{document}